%% file: DW-main.tex
\documentclass[pdflatex,sn-basic]{sn-jnl}

\usepackage{graphicx}%
\usepackage{multirow}%
\usepackage{amsmath,amssymb,amsfonts}%
\usepackage{amsthm}%
\usepackage{mathrsfs}%
\usepackage[title]{appendix}%
\usepackage{xcolor}%
\usepackage{textcomp}%
\usepackage{manyfoot}%
\usepackage{booktabs}%
\usepackage{algorithm}%
\usepackage{algorithmicx}%
\usepackage{algpseudocode}%
\usepackage{listings}%

\usepackage{todonotes}
\usepackage{tabu}
\usepackage{mathtools}
\usepackage{float}
\usepackage[normalem]{ulem}
\usepackage{xfrac}
\usepackage{subcaption}
\usepackage[table]{xcolor} 
\usepackage{mathtools, cuted}
\usepackage{placeins}

\theoremstyle{thmstyleone}%
\theoremstyle{thmstyletwo}%

\theoremstyle{thmstylethree}%

\begin{document}

\title[]{A Dantzig-Wolfe Reformulation for Automated Aircraft Arrival Routing and Scheduling}


\author*[1]{\fnm{Roghayeh} \sur{Hajizadeh}}\email{roghayeh.hajizadeh@liu.se}

\author[2]{\fnm{Tatiana} \sur{Polishchuk}}\email{tatiana.polishchuk@liu.se}

\author[1]{\fnm{Elina} \sur{R{\"o}nnberg}}\email{elina.ronnberg@liu.se}

\author[2]{\fnm{Christiane} \sur{Schmidt}}\email{christiane.schmidt@liu.se}

\affil[1]{\orgdiv{Department of Mathematics}, \orgname{Link{\"o}ping University}, \orgaddress{\city{Link{\"o}ping}, \postcode{581 83},  \country{Sweden}}}

\affil[2]{\orgdiv{Department of Science and Technology}, \orgname{Link{\"o}ping University}, \orgaddress{\city{Norrk{\"o}ping}, \postcode{601 74}, \country{Sweden}}}


\abstract{We consider the problem of computing dynamic arrival routes (DARs) in a terminal maneuvering area (TMA) together with an automated scheduling of all the arrivals within a given time interval.
	DARs are aircraft arrival routes from all entry points of a TMA, which merge traffic in that time interval towards the runway, 
	and are modeled as energy-efficient continuous-descent operations,  such that separation based on wake-turbulence categories is guaranteed within the TMA.  
	We propose a new model based on a Dantzig-Wolfe reformulation of a previous model for this problem.
	As in the previous model, we include tree consistency across consecutive planning intervals.
	However, the reformulation enables us to further improve the model and also consider aircraft that remain in the TMA from the previous period, a feature critical for operational safety. 
	In computational experiments for Stockholm Arlanda airport, the new model consistently outperforms the previous one: 
	we obtain solutions within 5 seconds to 12.65 minutes compared to up tp 40.9 hours with the old model for instances of half hours with high traffic.
	In addition,  we are able to solve instances of a full hour of arriving aircraft with high traffic (33 aircraft) within 22.22 to 58.57 minutes, whereas the old model could not solve these instances at all.
	While we schedule all aircraft as continuous-descent arrivals, our model can be applied to any type of speed profiles for the arriving aircraft.
}




\maketitle

\input{intro.tex}

\input{related-work.tex}

\input{prob-descrip.tex}
\input{model.tex}

\input{comput-res.tex}
\input{conclusion.tex}

 \backmatter

%
%
%

\bmhead{Acknowledgements}
This research was supported by grant 2022-03178 (OPT@ATM: Optimization Methods for Large Air Traffic Management Problems) from the Swedish Research Council (Vetenskapsrådet).
Some computational experiments were enabled by resources provided by the National Academic Infrastructure for Supercomputing in Sweden (NAISS), partially funded by the Swedish Research Council through grant agreement no. 2022-06725.

\section*{Declarations}


\begin{itemize}
	\item[] \textbf{Funding:} This research was supported by grant 2022-03178 (OPT@ATM: Optimization Methods for Large Air Traffic Management Problems) from the Swedish Research Council (Vetenskapsr{\aa}det).
	\item[] \textbf{Conflict of interests:} The authors have no relevant financial or non-financial interests to disclose.
    \item[] \textbf{Data availability:} The data will be made available upon request.
	\item[] \textbf{Author contribution:}\\
	\textit{Roghayeh Hajizadeh:} Data curation, Formal analysis, Investigation, Methodology, Software, Validation, Visualization, Writing--original draft, Writing--review and editing\\
	\textit{Tatiana Polishchuk:} Conceptualization, Data curation, Validation, Writing--review and editing\\
	\textit{Elina R{\"o}nnberg:} Conceptualization, Formal analysis, Methodology, Supervision, Writing--review and editing\\
	\textit{Christiane Schmidt:} Conceptualization, Formal analysis, Funding acquisition, Methodology, Project administration, Resources, Supervision, Writing--original draft, Writing--review and editing\\
\end{itemize}

\bibliography{TR-C-Dynamic-routes}

\begin{appendices}
	\input{appendixA.tex}
	
\end{appendices}

\end{document}

%% file: intro.tex
\section{Introduction}\label{intro-sec}
Air traffic volumes have increased for decades, 
in 2019, 9.2 billion passengers globally were reached according to~\cite{aci-gpter-23}. 
The numbers dropped significantly because of the Covid-19 pandemic,  but the International Air Transport Association~\cite{iata-gatdc-25} reported that air traffic reached 3.8\% above pre-pandemic levels in 2024,
and~\cite{iata-goat-23} predicts that air travel demand will double by 2040, growing at an average annual rate of 4.3\%.
This rise in air traffic volumes leads to an increased environmental impact and adds significant complexity to the tasks of air traffic controllers (ATCOs). 
These issues are of great importance  in terminal maneuvering areas (TMAs)---the airspace surrounding one or more airports---where all air traffic converges, leading to significant congestion and noise.
To manage the increasing volumes, it is essential to alleviate the environmental effects and the ATCO workload by providing improved arrival and departure procedures, which still enable a high runway throughput.

Continuous descent operations (CDOs) can reduce environmental impact. These operations, which involve optimal engine-idle descents, can lead to a decrease in fuel consumption, gaseous emissions, noise and fuel costs, see~\cite{e-ccdo-18}. 
Each aircraft has specific capabilities, and CDOs are adjusted accordingly, leading to unique optimal trajectories for each aircraft (where a trajectory gives both spatial and temporal position information, i.e., a 4D arrival).
However, these optimal trajectories often conflict with the strategic standard terminal arrival routes, 
resulting in decreased vertical and temporal predictability.  This situation prompts ATCOs to increase separation buffers, which negatively impacts throughput. Alternatively, ATCOs may issue instructions that modify the optimal trajectories, negatively affecting the environmental advantages.
In order to apply CDOs effectively, there is a need for tools that provide automated separation to ATCOs,  preferably  tools that provide an optimal flight sequence to achieve a high throughput (increasing the TMA's capacity).
As traffic demand varies over time, such tools must allow arrival routes to change during the day, leading to dynamic arrival routes (DARs).

For such a tool to be applicable in operations,  all possible factors influencing the situation must be taken into account (e.g., the aircraft type of an approaching aircraft), 
and the computation of the arrival routes and aircraft trajectories must be quick.  In a series of papers (\cite{atmos16,dgpss-cdadd-18, Saes_2019,Saez_2020,Polishchuk_2020,spshspp-asmda-21}), a team of authors  developed a mixed-integer-programming-based (MIP-based) model to compute optimal DARs with fully automated scheduling of CDOs that ensure aircraft separation and proposed an operational concept that facilitates the use of these routes.
For the computed trajectories, all aircraft fly according to their optimal neutral CDO speed profiles (which are generated beforehand), and an aircraft's arrival at the TMA entry point can be adjusted within a specific time window. 
The model selects the correct speed profile based on the length of the arrival route from the entry point to the runway. 
Furthermore, the progress of aircraft along routes within the given grid graph is monitored.
Despite demonstrating the feasibility of this approach, the framework suffers from long runtimes:
generating the arrival trees and sequence for a half-hour scenario for Stockholm Arlanda airport took between 1.58 hours for low-traffic cases and up to 40.9 hours for high-traffic cases, 
where light aircraft were artificially added to the traffic flow even though such aircraft are are hardly present at Arlanda at all.
This runtime is not practical for real-world operations, where DARs need to be recomputed approximately every 30 minutes, that is, within the time an aircraft typically spends in the TMA.
Hence, with that framework, it is also difficult to include additional features. For example,  \cite{Polishchuk_2020} and \cite{spshspp-asmda-21} did not integrate a critical aspect for operational safety---ensuring separation with aircraft that are still in the TMA but were scheduled in the previous time period.

Therefore, in this paper, we aim to improve the mathematical model to reduce runtimes. 
We again consider the speed profiles as input to our model. Any type of speed profile can be used, here, we use CDO speed profiles for all aircraft. 
We move from the arc- and flow-based model, where we assign  aircraft with a specific speed profile to a node in an overlayed grid at a specific point in time, to a path-based model, 
where we decide  on which path from each TMA entry point to the runway should be used---which then directly fixes the speed profile for each aircraft's trajectory.  
Hence, we introduce a Dantzig-Wolfe reformulation, based on discretization,  of the complete previous model 
and demonstrate that this can significantly reduce runtimes. 
Dantzig-Wolfe decomposition is a powerful technique used to reformulate a MIP problem from an original compact formulation to an extended formulation in a higher-dimensional space. 
The extended formulation is always at least as strong as the compact version and, in some cases, significantly stronger. 
This method has been widely applied in contexts such as vehicle routing, crew scheduling, and cutting stock problems
where it often leads to significantly reduced runtimes for challenging large-scale discrete optimization problems. 
For comprehensive introductions to this approach, see~\cite{JcLuGuBe2024} and~\cite{UcPeMo2024}.

This reformulation allows us to address more practical aspects (we integrated wind, \cite{hprs-dwraa-25-4}): a specific selected path again fixes a speed profile for an aircraft (which then depends on the wind direction w.r.t.~different segments on the path).
While we consider arrival paths in 2D, the aircraft trajectories (given by path plus CDO speed profile) determines both the spatial and temporal position of the aircraft in the TMA. Hence, they reflect the 4D CDO arrivals of all aircraft.

Our contributions within this paper are: 
\begin{itemize}
    \item We give a method to generate all possible paths from entry points to the runway.
    \item We present a path-based model, which is equivalent to the prior arc- and flow-based model. 
    \item We evaluate the model performance with real-world arrival data for Stockholm Arlanda airport. 
    \item We significantly improve runtimes (from hours to minutes) and solve instances that are not solvable with the previous model.  These runtimes are low enough for an operational context, however, operational usage requires the integration of additional factors, as we already did for wind (\cite{hprs-dwraa-25-4}).
    \item When we compute the arrivals for a certain time interval, some aircraft from the previous period are still present in the TMA. 
We incorporate this by either directly accounting for the presence of aircraft from the previous time or by enforcing a tree consistency between consecutive trees.
We then evaluate the performance of the model including these constraints on the same instances to illustrate the impact on the solutions.
\end{itemize}

{\bf Roadmap.} We introduce related work in Section~\ref{sec:rw} and
a problem definition in Section~\ref{prob-sec}. 
We present the previous model and our reformulated model in Section~\ref{model-sec}.
Furthermore, we give experimental results in Section~\ref{comp-sec}, and conclusions in Section~\ref{con-sec}. 
Finally, we provide additional tables in Appendix~\ref{appA}.

%% file: related-work.tex
\section{Related Work}\label{sec:rw}

In this section, we discuss related work in two directions that are relevant for our contribution. 
First, we review optimization approaches for ATM, with a focus on models that combine routing, sequencing, and trajectory decisions in a TMA. 
Second, we summarize methodological work on Dantzig-Wolfe reformulation and related column generation techniques, which provides the mathematical background for our path-based model.

\subsection{Related Work in Air Traffic Management}\label{subsec:rw-atm}

While there exist a relatively broad body of work on the optimization of aircraft trajectories (e.g., \cite{dp-ftsfc-15}), routing (e.g., \cite{atmos16}) and aircraft sequencing/separation (e.g., \cite{ishz-tpras-12}),  we focus on the integration of these---an area that has seen considerably less attention by researchers---and limit our review of related work to combined approaches. 

\cite{crmc-dors-10} presented a genetic-algorithm approach to compute aircraft arrival routes and the arrival sequence. First, they  developed distinct route topologies and then evaluated those with the heuristic-based scheduler.

\cite{krfs-roiat-23} presented an optimization framework that optimizes both aircraft trajectories and arrival sequence, while computing solutions that are robust against uncertainties from weather forecasts.  
In the framework, Kamo et al. employed robust optimal control and mixed-integer nonlinear programming (where they relaxed the integrality constraints to reduce runtimes). 
They aimed to develop functionality for a future arrival manager (AMAN) and reported that today's AMANs---an ATCO support tool---provide sequencing. Moreover,  \cite{e-aman-10} stated that a possible AMAN development (extended AMAN, E-AMAN) could lead to aircraft communicating a prefered trajectory, the estimated time of arrival (ETA), and a feasible time window around the ETA to ATC for optimization of the sequence. 
The pilots then use the calculated time of arrival sent by ATC as a required time of arrival for the flight management system to generate their trajectory. \cite{krfs-roiat-23} aimed for a robust optimization, as a previous deterministic optimization in their earlier work (\cite{krfs-ffp4d-22}) could lead to violation of operational constraints (like speed limitation) because of uncertainty. 
 They use two objective functions: maximizing the throughput and minimizing the fuel consumption.  They considered a very limited case study with only three approaching aircraft to the Leipzig/Halle airport. They employed six optimization runs (to avoid local optima), however, their relaxed approach using continuous variables yielded an optimal solution in only one out of six runs, while their approach using binary variables converged to the optimal solution in five out of six runs. 
 Moreover, the latter is less dependent on a (good) initial guess. For the very limited case study, both approaches have a runtime below a minute, the relaxed approach of about a second---however, the authors did not provide any larger, and thereby more realistic, traffic scenarios.

\cite{mdssa-iotma-19} integrated the operations at an airport and its surrounding terminal airspace: using a simulated-annealing heuristic, they decided on speed, arrival and departure times,  the runway assignment, and the pushback time at the airport.  To reduce computational time, Ma et al.~decomposed the problem using a sliding time window.  

Recently, \cite{ldczn-cppms-24} aimed to compare point merge and trombone arrival routes in TMAs and designed a mathematical model that decides on aircraft entry times and speeds and flight routes,  while making sure that no conflicts appear at merge points and on route segments.  The objective is to minimize delay and speed deviation for all arriving flights. 
Liu et al.  emphasized the importance of designing an efficient arrival route structure as route topologies in TMAs impact ``the scheduling of incoming aircraft flow and [can] increas(...)[e]  airspace capacity if properly addressed.'' 

\cite{hps-atstm-22} 
also focused on point merge: they presented a (simpler) MIP to sequence and merge the arriving traffic, where all arriving aircraft employ CDOs. Later, \cite{hps-aopmd-23} extended their MIP approach to also handle departing traffic in a dual-runway environment.

\cite{sgh-eocbs-24} provided an optimal-control-based model to estimate the fuel consumption of individual aircraft, which should be coupled with an approach to generate disjoint trajectories in a network by~\cite{hl-irhar-23} to find 4D TMA trajectories that balance fuel consumption and schedule deviation.
\cite{c-psmet-22} aimed to optimize fuel consumption for arriving aircraft using path stretching to sequence the arriving traffic.  He provided a MIP model and compared the performance to point merge for Istanbul Sabiha Gök\c{c}en airport.
\cite{t-amoam-19} aimed to optimize both aircraft sequencing and trajectories and formulated a MIP model and an optimal control problem, respectively.
\cite{sddp-oasrt-14} considered aircraft scheduling and routing using a MIP model. However, they do not take aircraft-specific performance into account.
The combined assignment of routes and times has been applied in non-aviation context, e.g.,  in dynamic vehicle routing by~\cite{br-sdvrp-91}.

\subsection{Dantzig--Wolfe Reformulation and Related Methods}\label{subsec:rw-dw}
Dantzig-Wolfe decomposition is a technique used to reformulate a MIP problem from an original compact formulation
to an extended formulation in a higher-dimensional space. 
The extended formulation will always be at least as strong as the compact formulation, and sometimes it is much stronger. 
In such cases the reformulation can contribute to significantly reduced runtimes for challenging large-scale discrete optimisation problems.
In this section, we explain the technical procedure of Dantzig-Wolfe reformulation for a special case of MIP problems: those with binary variables and a block diagonal structure.

Consider the compact binary linear programming problem

{\small
\begin{align}\label{dw-1}
	\min~~ z =& \sum_{k \in K} c_k^Tx^k\\\label{dw-2}
	\text{s.t.}~~~~~~~&\sum_{k \in K}A^{(1)}_kx^k\geq b^{(1)},\\\label{dw-3}
	&A^{(2)}_kx^k\geq b^{(2)}_k,&&\forall k\in K,\\\label{dw-4}
	&x^k~~\text{binary},&&\forall k\in K,%
\end{align}
}%
where $K$ is a set of indices and  $x^k$, $k\in K$, is a vector  consisting of $n^k$ variables. 
Moreover, $c_k$, $A^{(1)}$, $A^{(2)}_k$, $b^{(1)}$, and $b^{(2)}_k$  are of compatible size and $\mathbf{x}=(x^k)_{k\in K}$ is the vector of variables that includes all $x^k$. 
This model has a block diagonal structure where all the variables are connected in~\eqref{dw-2} and  separated over $k$ in~\eqref{dw-3}.

For each $k\in K$, there exists a finite set of binary points 
$\{\bar{x}_q^k\}_{q\in Q^k}$ that enumerates the solutions in the set $\mathcal{X}^k=\{x^k: A^{(2)}_kx^k\geq b^{(2)}_k, ~ x^k ~ \mbox{binary} \}$.
Thus, any $x^k\in\mathcal{X}^k$ can be written as
{\small
\begin{align}\label{dw-5}
~~~~~~x^k=&\sum_{q\in Q^k}\bar{\lambda}_q^k\bar{x}_q^k,\\\label{dw-6}
	&\sum_{q\in Q^k}\bar{\lambda}_q^k=1,\\\label{dw-7}
	&\bar{\lambda}_q^k~~\text{binary},&&\forall q\in Q^k.
\end{align}}
According to Constraints~\eqref{dw-6}--\eqref{dw-7}, $x^k$ in Constraint~\eqref{dw-5} is equal to one of the points from $\{\bar{x}_q^k\}_{q\in Q^k}$.
Note that these solutions together give the feasible solutions to constraints~\eqref{dw-3}--\eqref{dw-4}. 
By substituting $x^k$ with~\eqref{dw-5}--\eqref{dw-7} in the compact  model~\eqref{dw-1}--\eqref{dw-4}, we obtain  an extended reformulation 

{\small
\begin{align}\label{dw-8}
	\min~~ z =& \sum_{k \in K} \sum_{q\in Q^k}\bar{c}_q^k\bar{\lambda}_q^k \\\label{dw-9}
	\text{s.t.}~~~~~~~& \sum_{k \in K} \sum_{q\in Q^k}\bar{a}_{q}^k\bar{\lambda}_q^k\geq b^{(1)},\\\label{dw-10}
	&\sum_{q\in Q^k}\bar{\lambda}_q^k=1,&&\hspace*{-1.5cm}\forall k\in K,\\\label{dw-11}
	&\bar{\lambda}_q^k ~~\text{binary},&&\hspace*{-1.5cm}\forall k\in K, ~\forall q\in Q^k,
\end{align}
}%
where $\bar{c}_q^k=c_k^T\bar{x}_q^k$ and $\bar{a}_{q}^k= A_k^{(1)}\bar{x}_q^k$. 
In the compact model,  the $x^k$ are the variables while in the extended formulation, the ${\lambda}_q^k$ are the variables. 
When the reformulation is based on a complete enumeration of solutions in $\mathcal{X}^k$, the approach is called Dantzig-Wolfe reformulation based on discretization. 
Another procedure for Dantzig-Wolfe reformulation is convexification, where instead of enumerating all solutions, extreme points and extreme rays are enumerated.
When the variables are binary, convexification and discretization coincide, as remarked by~\citet[p.~1012]{LuDe05}.

The extended formulation is at least as strong as the original model according to~\cite{Vand00}, i.e.,\ 
the bound obtained by the linear programming relaxation of the extended model is at least as tight as the bound obtained by the linear programming relaxation of the original model. 
However, in practice, it can be much stronger for certain types of problem structures.

In optimization problems like ours, which can be represented as graphs, the original model often makes decisions based on arcs.
This means that the variables correspond to the presence or absence of certain arcs in the solution 
and the constraints ensure that the selected arcs form a valid solution according to the problem's requirements such as forming a path. 
However, in the extended formulation, paths are constructed and decisions are made based on these paths. 
The variables in this model correspond to selecting these paths rather than individual arcs.
This shift from arcs to paths is characteristic of Dantzig-Wolfe approaches and appears in classical applications such as vehicle routing, crew scheduling, and multicommodity flow, see~\cite{JcLuGuBe2024} for more.

It is very common to apply column generation and branch-and-price to these types of 
reformulated models since the number of variables in the extended formulation typically grows exponentially.
For our problem,  each variable of the extended formulation corresponds to a path.
We show in Subsection~\ref{subsec:gen-paths} that 
it is possible to enumerate all the paths for a discretization of the TMA into a relatively coarse grid. 
Therefore, it is possible to construct the full path-based extended formulation without applying column generation or branch-and-price.
The feasibility of explicit path enumeration in our model is tied to the current level of detail in the model.
If the grid becomes finer, the number of feasible paths will grow rapidly.
In that case one would no longer be able to enumerate all paths, and the natural next step would be to develop a full branch-and-price-and-cut framework,
with a pricing problem generating promising trajectories on demand.
Thus, the present formulation can be seen as a step toward more advanced algorithmic developments.

%% file: prob-descrip.tex
\section{Problem Description}\label{prob-sec}
Given a set of aircraft entering a TMA through specified entry points aiming to land on a runway, 
the problem is to determine an optimal arrival tree for these aircraft  aiming to merge traffic seamlessly from the entry points to the runway.
 The input includes location of TMA entry points and the runway for an airport,  and a set of aircraft planned to arrive within a certain time window to their entry points.
 Each entry point provides multiple possible paths with varying distances to the runway. These paths merge at intermediate points within the TMA before meeting at the final merging fix (a precisely defined navigation point in the sky which is used to guide and sequence aircraft), from where they will proceed to the runway.
 Since merge points require more ATCO attention, the lowest possible traffic complexity is required, \cite{pkmkz-fprps-08}. 
 High traffic complexity can be reduced by limiting the number of routes merging at a point and by restricting the number of merging points within a small area.
 An arrival tree, directed from the entry points (leaves) to the runway (root), need to comply with these operational constraints:

\begin{enumerate}
	\item No more than two routes merge at a point. 
	\item  Merge points are separated by a minimum distance $L$, \cite{pkmkz-fprps-08}. 
	\item Routes do not make sharp turns in order to respect aircraft dynamics, that is,  any turn from a route segment to the consecutive route segment must be larger than a given minimum  angle $\gamma$.
	\item Obstacles, like no-fly zones and noise-sensitive areas, are avoided.
	\item To ensure safety,  all aircraft are temporally separated along the arrival routes. 
	The required temporal separation between a leading and trailing aircraft depends on their respective wake-turbulence categories.
	\item All aircraft arrive to their entry points within a given interval  around the planned time. 
	That is, in order to enable the accommodation of more aircraft within the scheduling framework,  deviations from the planned arrival times at the entry points are allowed. 
	\item All aircraft follow a set of  speed profiles that depend on the arrival route length. In our framework, they follow CDO speed profiles. 
\end{enumerate}

The goal is to design DARs that balance efficiency and environmental considerations. 
One objective is to minimize the total distance covered by all aircraft in the TMA, ensuring shorter flight routes for aircraft. 
Additionally, the arrival route tree should occupy a minimal area to reduce ATCOs' workload and limit the spread of noise and other environmental impacts. 
These objectives, referred to as paths length and tree weight, respectively, are both considered in the mathematical model. 

%% file: model.tex
\section{Mathematical Formulation}\label{model-sec}
This section introduces the Dantzig-Wolfe reformulation of the compact model originally proposed by~\cite{spshspp-asmda-21}. 
The original model is arc- and flow-based while our reformulation yields a path-based model, offering new insights and computational advantages. 
To provide a comprehensive understanding, in Subsection~\ref{subsec:grid-turbulence}, we outline the general problem specifics, 
including a detailed explanation of the grid discretization of the TMA, arrival-time constraints, and wake-turbulence categories.
This is the same grid discretization as the one used by~\cite{spshspp-asmda-21}. 
Following this, in Subsection~\ref{subsec:old-model}, we give an overview of the compact model proposed by~\cite{spshspp-asmda-21}, providing a foundation for understanding the subsequent reformulation. 
Finally, we present the path-based formulation of the problem in Subsection~\ref{subsec:new-model}
and discuss consistency between trees of different time periods in Subsection~\ref{subsec:consistency}.

\subsection{Grid Overlay, Arrival Times, and Wake-Turbulence Categories}\label{subsec:grid-turbulence}
\begin{figure} 
	\begin{subfigure}[t]{0.5\textwidth}
	\centering
	\includegraphics[width=0.74\textwidth]{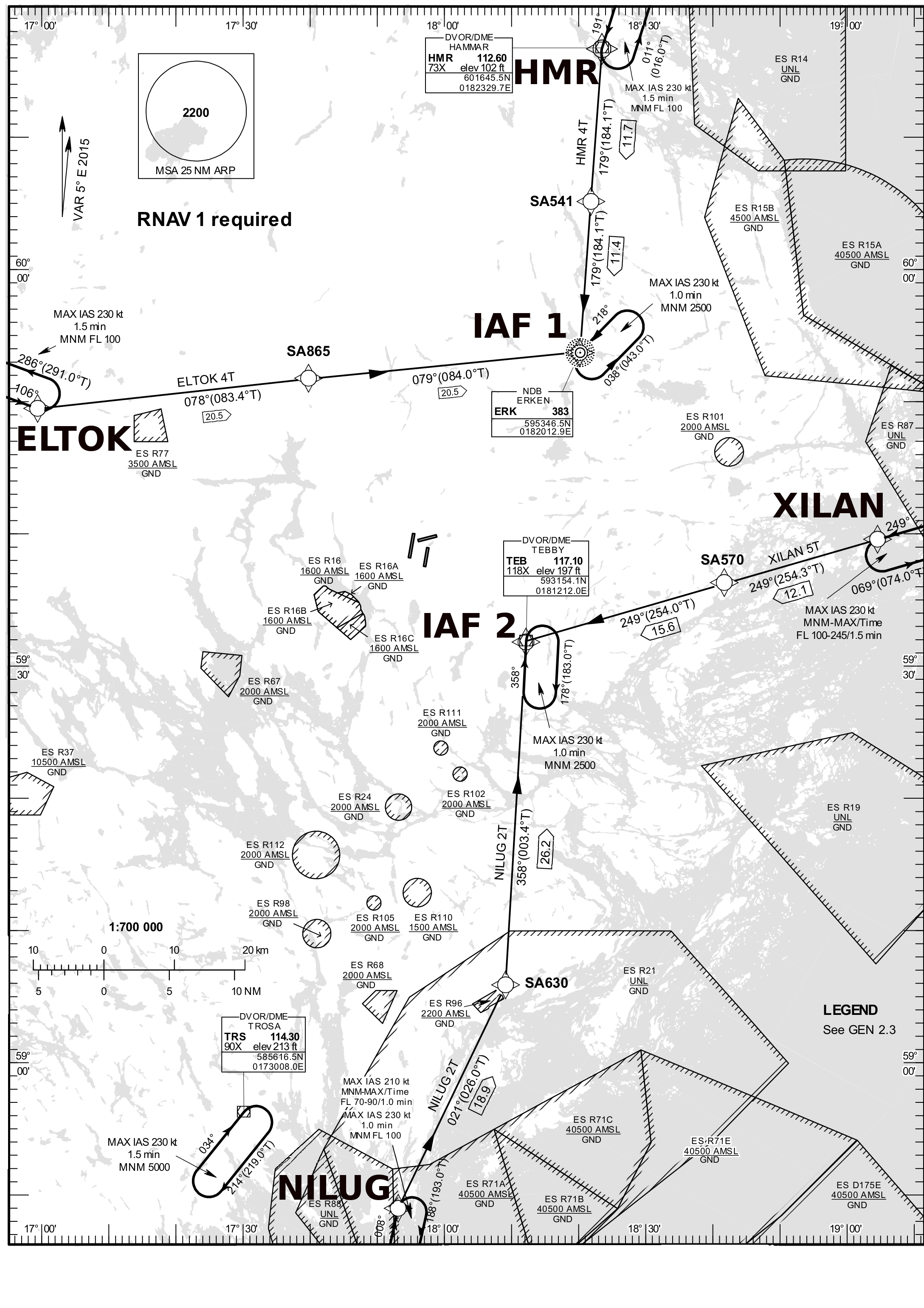} \hspace*{1cm}\subcaption{Arlanda TMA with runway 19L \\ and four entry points}
\end{subfigure}\hfil\hspace*{-1.5cm}
\begin{subfigure}[t]{0.5\textwidth}
	\centering
	\includegraphics[width=1\textwidth]{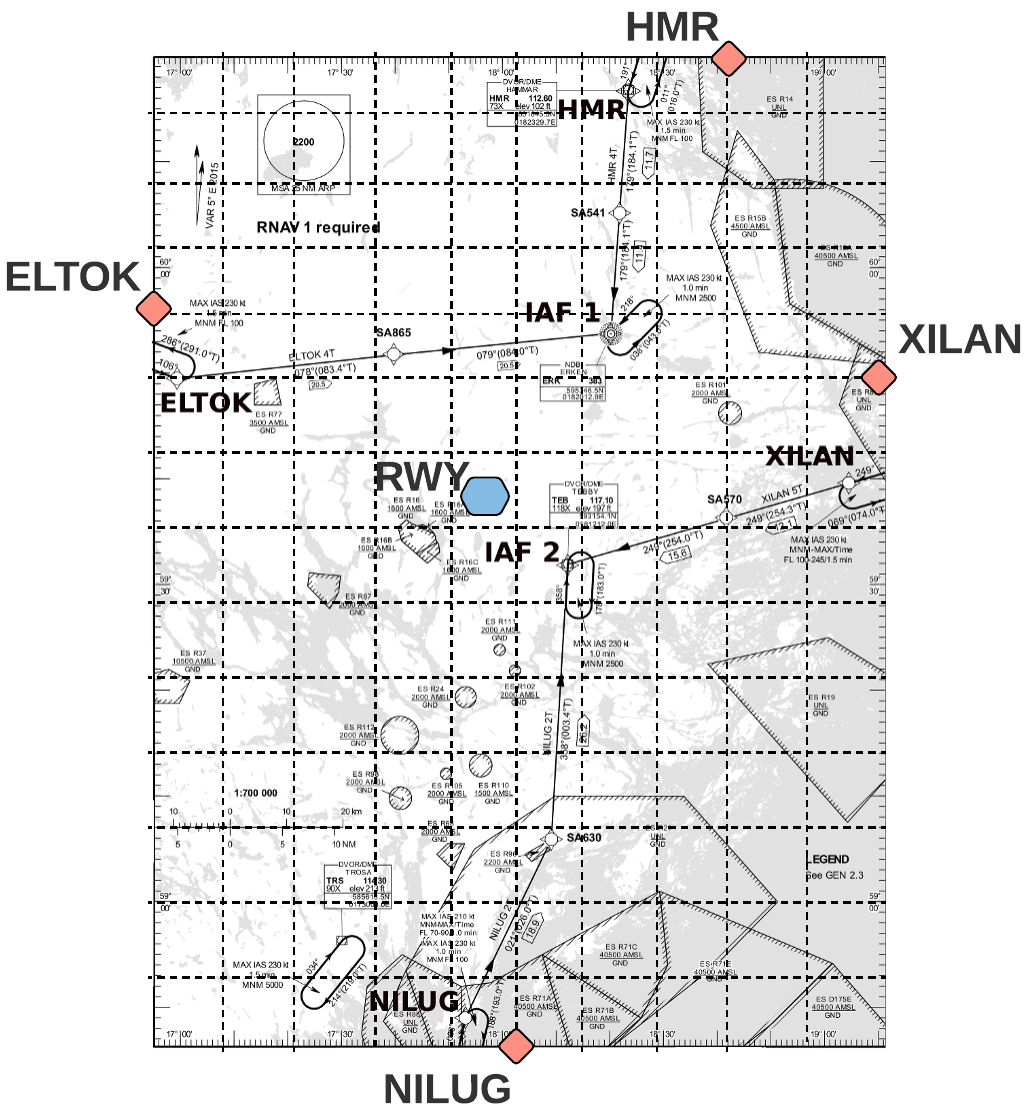} \hspace*{.7cm}\subcaption{$15\times 11$ grid over Arlanda TMA}
\end{subfigure}\hspace*{-.75cm}
\caption{Overlaying Arlanda TMA with a $15\times 11$ grid}\label{fig:arlanda-TMA-grid} 
\end{figure}

We discretize the TMA into a square grid of size $q\times n$, 
where the side length of each grid pixel is equal to the minimum separation distance $L$.  This  pixel side length also yields the desired merge-point separation.
The positions of both the entry points and the runway are then adjusted to align with this grid. 
For instance, the TMA of Stockholm Arlanda airport given in Figure~\ref{fig:arlanda-TMA-grid}(a) is overlaid with a $15\times 11$ grid in Figure~\ref{fig:arlanda-TMA-grid}(b) which ensures a separation of about 6NM.

From this grid, we derive a bi-directed graph  $G=(V,E)$ with nodes $V$ and edges $E$, where each grid point corresponds to a node in $V$.
In the graph, every node is connected to its eight neighbors,
and for any two neighboring nodes $i$ and $j$, both edges $(i,j)$  and $(j,i)$ are included in $E$.
Entry points and the runway are exceptions: entry points do not have incoming edges and the runway does not have outgoing edges.
We denote the length of edge $(i,j)\in E$  by $l_{ij}$, the set of entry points by $\mathcal{P}\subseteq V$, and the runway by $r\in V$.

In addition, we denote the set of all aircraft arriving at entry point $b\in\mathcal{P}$ by $\mathcal{A}_b$,  
and the set of all aircraft by $\mathcal{A} =\bigcup_{b\in \mathcal{P}}\mathcal{A}_b$. 
Furthermore, let $\Gamma_{ij}$ for each edge $e=(i, j)$ be the set of all outgoing edges from node $j$ that form an angle smaller that $\gamma$ with $e$.

Moreover, assume that  $\bar{t}_a$ denotes the planned arrival time of aircraft $a$ at its entry point.
To ensure the feasibility of the required temporal separation, we allow deviations from the planned
time at the entry points so that aircraft $a$ can arrive sometime in the interval  $[\bar{t}_a-\mu,\bar{t}_a+\mu]$, 
where $\mu$ denotes a small integer, that is, denoting a number of minutes. 
Moreover, we let $\overline{T}$ be an upper bound on the time all considered aircraft spend in the TMA until landing.
We discretize the time into periods from $0$ to $\overline{T}$ and denote the set of time intervals by $T=\{0,\ldots,\overline{T}\}$.

As mentioned in Section~\ref{prob-sec}, the temporal separation requirement depends on the wake-turbulence categories of the leading and trailing aircraft.
Assume that $\mathcal{C}=\{1, 2, \ldots, n_c\}$ is the set of indices of all different disjoint wake-turbulence categories 
where $C_{\kappa}$  for $\kappa\in\mathcal{C}$ is a set of aircraft from the same category and $\mathcal{A}=\bigcup_{\kappa\in \mathcal{C}} C_{\kappa}$.
Let $\sigma_{\kappa_1, \kappa_2}$ denote the required temporal separation if the leading aircraft is from $C_{\kappa_1}$ 
and the trailing aircraft is from $C_{\kappa_2}$.\
We let $\Omega=\max_{\kappa_1,\kappa_2\in \mathcal{C}}~~ \sigma_{\kappa_1, \kappa_2}$.

Let $\lambda$ denote the upper bound on the number of nodes in any path, and define $\mathcal{L}=\{1,\ldots,\lambda\}$ as the set of all possible numbers of nodes that a path can have.
In our context, a route is a sequence of nodes connected by edges. To emphasize the order in which these connections are traversed, we refer to each edge in a route as a segment. 
For aircraft $a$,  $\mathcal{S}(a)$ represents the set of CDO speed profiles,  with each profile corresponding to a route of a particular length. 
Each speed profile specifies the time required to traverse the first, second, and subsequent segments of the route, where the time to cover a segment depends on the number of segments completed along the route up to the current node. 
We note that wind effects are not modeled here, so profiles are assumed to depend only on route length and position, although in practice wind can cause profiles for routes of the same length to differ.
In particular, let  $u_{a,p,k}$ be the time that aircraft $a$ using speed profile $p$ needs to traverse segment number $k$ on its route from entry point to the runway.

\subsection{Arc- and Flow-Based Formulation}\label{subsec:old-model}
In this section, we describe the mathematical model introduced by \cite{spshspp-asmda-21}, which forms the basis for optimizing DARs planning in the TMA.  
To achieve this, the following variables are introduced.
\begin{itemize}
	\item Binary variables $x_{i,j}$  indicate whether  or not edge $(i,j)$ is included in the arrival tree.
	\item Flow variables $f_{i,j}$ give the flow on edge $(i,j)$.
	\item Binary variables ${x}^{\scalebox{0.7}{EP}}_{(i,j),b}$ indicate whether or not edge $(i,j)$ participates in the route from entry point $b$ to the runway.
	\item Variable $\ell(b)$ gives the length of the path from entry point $b$ to the runway in the arrival tree, measured by the number of grid edges rather than the physical distance.
	\item Binary variables  $y_{a,j,p,k,t}$ indicate whether or not  aircraft $a$ using speed profile $p$ occupies the $k$-th node $j$ at time $t$.
	 \item Binary variable $z_{a,j,i,b,p,k,t}$ takes value 1 if and only if aircraft $a$ using speed profile $p$ occupies node $j$ at time $t$ as the $k$-th node in its route, and the edge $(j,i)$ is part of the  path from entry point $b$ to the runway. Otherwise, $z_{a,j,i,b,p,k,t}=0$.
	\item Binary variables $\psi_{b,a,p}$ and $\phi_{b,a,p}$ control whether the length of the speed profile $p$ for aircraft $a$ is equal to $\ell(b)$.
\end{itemize}
We call this model M1 and it is given as:

\setcounter{equation}{-1}  
\makeatletter
\renewcommand{\theequation}{M1.\arabic{equation}}
\makeatother

{\allowdisplaybreaks
	\footnotesize
		\centering
		\begin{align}
			\min \quad     		 & \beta \sum_{(i,j)\in E} \ell_{i,j} x_{i,j} + (1-\beta)\sum_{(i,j)\in E} \ell_{i,j} f_{i,j}   &&\label{obj-1}\\
			\mbox{s.t.}\quad& \sum_{k:(k,i)\in E} f_{k,i} - \sum_{j:(i,j)\in E} f_{i,j}  = &&\hspace*{-3.75cm} \begin{cases} |\mathcal{A}|, \;\;\;\;\;\;\;\;\;\;\;  i=r, \\ -|\mathcal{A}_i|,\;\;\;\;\;\;\;   i\in\mathcal{P},\\ 0,  \;\;\;\;\;\;\;\;\;\;\;\;\;\; i\in V\setminus\{\mathcal{P}\cup r\}, \end{cases} \label{eq:flow-weighted}\\
			& x_{i,j} \geq \quad \frac{f_{i,j}}{Q}, &&\hspace*{-1.5cm}\forall (i,j)\in E, \label{eq:flow-u}\\
			& \sum_{k:(k,i)\in E} x_{k,i} \leq  2, &&\hspace*{-1.5cm} \forall i \in V\setminus\{\mathcal{P}\cup r\},\label{eq:mcf-indegree}\\
			& \sum_{j:(i,j)\in E} x_{i,j} \leq  1, &&\hspace*{-1.5cm}  \forall i \in V\setminus\{\mathcal{P}\cup r\},\label{eq:mcf-outdegree}\\
			& \sum_{k:(k,r)\in E} x_{k,r} = 1,&&  \label{eq:mcf-indegree-vstar}\\
			& \sum_{j:(i,j)\in E} x_{i,j} = 1, &&\hspace*{-1.5cm}  \forall i\in \mathcal{P},\label{eq:mcf-outdegree-S}\\
			& |\Gamma_{i,j}| x_{i,j} + \sum_{(i', j')\in \Gamma_{i,j}} x_{i',j'} \leq |\Gamma_{i,j}|, &&\hspace*{-1.5cm} \forall (i,j)\in E, \label{eq:curvature}\\
			&x_{i, i+1+n} + x_{i+1+n, i} + x_{i+n, i+1} + x_{i+1, i+n} \leq 1,&&\nonumber\\
			&&&\hspace*{-3cm}  \forall i \in V'\setminus\{\mathcal{P}\cup r\}: i+1+n, i+n, i+1 \not\in\{\mathcal{P}\cup r\}, \label{eq:aux-1}\\
			&x_{i, i+1+n} + x_{i+n, i+1} + x_{i+1, i+n} \leq  1, &&\hspace*{-1.5cm} \forall i \in \mathcal{P},    \label{eq:aux-2} \\
			&x_{i, i+1+n} + x_{i+1+n, i} + x_{i+1, i+n} \leq 1, &&\hspace*{-1.5cm}\forall i: i+1 \in \mathcal{P},     \label{eq:aux-3}\\
			&x_{i, i+1+n} + x_{i+n+1, i} + x_{i+n, i+1} \leq  1, &&\hspace*{-1.5cm}\forall i: i+n \in \mathcal{P},   \label{eq:aux-4} \\
			&x_{i+1+n, i} + x_{i+n, i+1} + x_{i+1, i+n} \leq  1, &&\hspace*{-1.5cm}\forall i: i+n+1 \in \mathcal{P}, \label{eq:aux-5}\\
			&x^{\scalebox{0.7}{EP}}_{(i,j),b} \leq\;  x_{i,j}, &&\hspace*{-1.5cm}\forall b\in\mathcal{P},\forall (i,j)\in E,\label{eq:xbs1}\\
			& \sum_{j:(b,j)\in E} {x}^{\scalebox{0.7}{EP}}_{(b,j),b}  =  1, &&\hspace*{-1.5cm}\forall b\in\mathcal{P}, \label{eq:xbs2}\\
			&\sum_{j:(j,r)\in E} x^{\scalebox{0.7}{EP}}_{(j,r),b}\;  =  1, &&\hspace*{-1.5cm}\forall b\in\mathcal{P},\label{eq:xbs3}\\
			&\sum_{i:(i,j)\in E} x^{\scalebox{0.7}{EP}}_{(i,j),b} -\sum_{k:(j,k)\in E} x^{\scalebox{0.7}{EP}}_{(j,k),b}  = 0, &&\hspace*{-1.5cm}\forall j\in V\setminus\{\mathcal{P}\cup r\}, \forall b\in\mathcal{P}, \label{eq:xbs4}\\
			&\sum_{t=\bar{t}_a-\mu}^{\bar{t}_a+\mu}\sum_{p\in\mathcal{S}(a)} y_{a,b,p,1,t}  = 1, &&\hspace*{-1.5cm} \forall b \in\mathcal{P}, \forall a\in \mathcal{A}_b, \label{eq:y1n_2}\\
			&y_{a,b,p,k,t}  =  0, &&\hspace*{-3.5cm} \forall b \in\mathcal{P}, \forall a\in\mathcal{A}_b, \forall p\in\mathcal{S}(a), \forall t\leq \bar{t}_a-\mu,\forall k\neq 1\in\mathcal{L}, \label{eq:y2n_2}\\
			& y_{a,b,p,1,\bar{t}_a}\;  = 0,&&\hspace*{-3.5cm} \forall b \in\mathcal{P}, \forall a\in \mathcal{A}_b,\forall p\in\mathcal{S}(a), t\in T: t\notin [\bar{t}_a-\mu, \bar{t}_a+\mu],\label{eq:y3n_2}\\ 
			&y_{a,b',p,k,t}  =  0, &&\hspace*{-3.5cm}  \forall b'\neq b \in\mathcal{P}, \forall a\in \mathcal{A}_b,\forall p\in\mathcal{S}(a), \forall k\in\mathcal{L},\forall t\in T, \label{eq:y4n} \\
			&y_{a',b,p,1,\bar{t}_a}  =  0, &&\hspace*{-3.5cm}  \forall b \in\mathcal{P}, \forall a'\neq a\in \mathcal{A}_b,\forall p\in\mathcal{S}(a),  \label{eq:y5n}\\
			&y_{a,j,p,1,t}  =  0, &&\hspace*{-3.5cm} \forall j \in V\setminus\mathcal{P}, \forall a\in \mathcal{A},\forall p\in\mathcal{S}(a), \forall t\in T, \label{eq:y7n} \\
			&y_{a,j,p,k,t}  =  0,&&\hspace*{-3.5cm}  \forall b \in\mathcal{P}, \forall a\in\mathcal{A}_b, \forall p\in\mathcal{S}(a), \forall t< \bar{t}_a-\mu, \forall k\in\mathcal{L}, \label{eq:nowhere-before}\\
			&y_{a,j,p,k,t}  \leq  \sum_{\substack{i\in V: (i,j)\in E}} x_{i,j}, &&\hspace*{-3.5cm} \forall j \in V\setminus\mathcal{P}, \forall a\in\mathcal{A}, \forall p\in\mathcal{S}(a), \forall k\in\mathcal{L}, \forall t\in T, \label{eq:y6n} \\
			&\sum_{p\in\mathcal{S}(a)}y_{a,j,p,k,t}  \leq  1, &&\hspace*{-3.5cm}  \forall j \in V, \forall a\in \mathcal{A},\forall k\in\mathcal{L},\forall t\in T, \label{eq:y8n}\\
			&\sum_{t=\bar{t}_a-\mu}^{\bar{t}_a+\mu} y_{a,b,p,1,t} =  1-\psi_{b,a,p}, &&\hspace*{-3.5cm}  \forall b \in\mathcal{P}, \forall a\in\mathcal{A}_b, \forall p\in\mathcal{S}(a),\label{eq:length-profile}\\
			&\ell(b) =  \sum_{(i,j)\in E} x^{\scalebox{0.7}{EP}}_{(i,j),b}, &&\hspace*{-1.5cm}\forall b\in\mathcal{P}, \label{eq:length}\\
			&-\psi_{b,a,p} \leq  (\ell(b)-p+\frac{1}{2})/\lambda +\lambda\cdot\phi_{b,a,p}, &&\hspace*{-1.5cm} \forall b \in\mathcal{P}, \forall a\in \mathcal{A}_b, \forall p\in\mathcal{S}(a),\label{eq:psi1}\\
			&\psi_{b,a,p}  \leq  -\ell(b)+p+\frac{1}{2}+\lambda\cdot\phi_{b,a,p},&&\hspace*{-1.5cm}\forall b \in\mathcal{P}, \forall a\in\mathcal{A}_b, \forall p\in\mathcal{S}(a),\label{eq:psi2}\\
			&\psi_{b,a,p}  \leq  (\ell(b)-p)+\frac{1}{2}+\lambda(1-\phi_{b,a,p}),&&\hspace*{-1.5cm}\forall b \in\mathcal{P}, \forall a\in \mathcal{A}_b, \forall p\in\mathcal{S}(a), \label{eq:psi3}\\
			&-\psi_{b,a,p}  \leq  (-\ell(b)+p+\frac{1}{2})/\lambda +\lambda(1-\phi_{b,a,p}),&&\hspace*{-1.5cm}\forall b'\neq b \in\mathcal{P}, \forall a\in \mathcal{A}_b,\forall p\in\mathcal{S}(a),\label{eq:psi4} \\[.5em]
			&\sum_{a\in B} \sum_{p\in\mathcal{S}(a)}\sum_{k\in\mathcal{L}}\sum_{\tau=t}^{t+\sigma_{AB}-1}   y_{a,j,p,k,\tau} \leq &&\hspace*{-3cm} \Omega - \Omega\cdot\sum_{a'\in A}\sum_{\substack{p'\in\mathcal{S}(a')}} \sum_{k'\in\mathcal{L}} y_{a',j,p',k',t}, &&\nonumber\\
			&&&\hspace*{-1.5cm}\forall j\in V, \forall t\in\{0,\ldots,\overline{T}-\sigma_{\kappa_1, \kappa_2}+1\}, \label{eq:wake-turb-1}\\[1em]
			&\sum_{a\in B} \sum_{p\in\mathcal{S}(a)}\sum_{k\in\mathcal{L}}\sum_{\tau=t}^{\overline{T}}   y_{a,j,p,k,\tau} \leq  &&\hspace*{-3.5cm}\Omega - \Omega\cdot\sum_{a'\in A}\sum_{\substack{p'\in\mathcal{S}(a')}} \sum_{k'\in\mathcal{L}} y_{a',j,p',k',t}, &&\nonumber\\
			&&&\hspace*{-1.5cm}\forall j\in V,\forall t\in\{\overline{T}-\sigma_{\kappa_1, \kappa_2}+2,\ldots,\overline{T}\}, \label{eq:wake-turb-1'}\\ 
			&\sum_{a\neq a'\in A} \sum_{p\in\mathcal{S}(a)}\sum_{k\in\mathcal{L}}\sum_{\tau=t}^{t+\sigma_{AA}-1}   y_{a,j,p,k,\tau} \leq  \Omega - \Omega\cdot y_{a',j,p',k',t},&&\nonumber\\
			&&&\hspace*{-3.8cm}\forall a' \in A, \forall p'\in\mathcal{S}(a'), \forall j\in V,\forall t\in\{0,\ldots,\overline{T}-\sigma_{\kappa, \kappa}+1\}, \label{eq:wake-turb-2}\\
			&\sum_{a\neq a'\in A} \sum_{p\in\mathcal{S}(a)}\sum_{k\in\mathcal{L}}\sum_{\tau=t}^{\overline{T}}   y_{a,j,p,k,\tau} \leq  \Omega - \Omega\cdot y_{a',j,p',k',t},&&\nonumber\\
			&&&\hspace*{-3.8cm}\forall a' \in A, \forall p'\in\mathcal{S}(a'), \forall j\in V,\forall t\in\{\overline{T}-\sigma_{\kappa, \kappa}+2,\ldots,\overline{T}\}, \label{eq:wake-turb-2'}\\
			& z_{a,j,i,b,p,k,t-u_{a,p,k}} \leq  x^{\scalebox{0.7}{EP}}_{(j,i),b},&&\nonumber\\
			&&&\hspace*{-5cm}\forall b\in\mathcal{P},\forall a\in\mathcal{A}_b, \forall (j,i) \in E,  \forall k\in\mathcal{L}, \forall t\in \{u_{a,p,k}+1,\ldots, \overline{T}+1\}, \label{eq:z1n}\\
			& z_{a,j,i,b,p,k,t-u_{a,p,k}} \leq   y_{a,j,p,k,t-u_{a,p,k}}, &&\nonumber\\
			& &&\hspace*{-5cm}\forall b\in\mathcal{P},\forall a\in\mathcal{A}_b, \forall (j,i) \in E,  \forall k\in\mathcal{L}, \forall t\in \{u_{a,p,k}+1,\ldots, \overline{T}+1\}, \label{eq:z2n}\\
			& z_{a,j,i,b,p,k,t-u_{a,p,k}} \geq   x_{(j,k),b} - (1 - y_{a,j,p,k,t-u_{a,p,k}}), &&\nonumber\\
			& &&\hspace*{-5cm}\forall b\in\mathcal{P},\forall a\in\mathcal{A}_b, \forall (j,i) \in E,  \forall k\in\mathcal{L},\forall t\in \{u_{a,p,k}+1,\ldots, \overline{T}+1\},  \label{eq:z3n}\\
			&\sum_{j: (j,i)\in E} z_{a,j,i,b,k,t-u_{a,p,k}} - y_{a,i,p,k+1,t} = 0,  &&  \nonumber\\
			& &&\hspace*{-5.8cm}\forall b\in\mathcal{P}, \forall a\in\mathcal{A}_b,\forall p\in\mathcal{S}(a), \forall k\in\mathcal{L},\forall i\in V\setminus\mathcal{P}, \forall t\in\{u_{a,p,k}+1,\ldots,\overline{T}\},\label{eq:yaktu-new}\\
			&x_{i,j}\quad\quad\quad~~~\quad\text{binary},&&\hspace*{-3.5cm}\forall (i,j)\in E,\label{eq:dom1}\\
			&x^{\scalebox{0.7}{EP}}_{(i,j),b}\quad~~\quad\quad\text{binary},&&\hspace*{-3.5cm}\forall (i,j)\in E, \forall b\in\mathcal{P},\label{eq:dom2}\\
			&y_{a,j,p,k,t}\quad\quad~~\text{binary},&&\hspace*{-3.5cm}\forall a\in\mathcal{A}, \forall p\in\mathcal{S}(a), \forall j \in V, \forall k\in\mathcal{L}, \forall t\in T,\label{eq:dom3}\\
			&z_{a,j,i,b,p,k,t}\quad~\text{binary},&&\hspace*{-3.5cm}\forall b\in\mathcal{P}, \forall a\in\mathcal{A}_b, \forall p\in\mathcal{S}(a), \forall (j,i) \in E, \forall k\in\mathcal{L}, \forall t\in T,\label{eq:dom4}\\
			&\psi_{b,a,p}\quad\quad\quad\quad\text{binary},&&\hspace*{-3.5cm}\forall b\in\mathcal{P}, \forall a\in \mathcal{A}_b, \forall p\in\mathcal{S}(a),\label{eq:dom5}\\
			&\phi_{b,a,p} \quad\quad\quad\quad\text{binary},&&\hspace*{-3.5cm}\forall b\in\mathcal{P}, \forall a\in \mathcal{A}_b, \forall p\in\mathcal{S}(a),\label{eq:dom6}\\
			&f_{i,j}\geq0,&&\hspace*{-3.5cm}\forall (i,j)\in E,\label{eq:dom7}\\
			&\ell(b)\geq0,&&\hspace*{-3.5cm}\forall b\in\mathcal{P}.\label{eq:dom8}
	\end{align}
		 }		 	

The objective function~\eqref{obj-1} is a convex combination of the weight of the arrival-route tree  and the total aircraft distances with coefficient $\beta$. 
Constraint~\eqref{eq:flow-weighted} is the node-equilibrium constraint which ensures that a flow of $|\mathcal{A}|$ reaches the runway $r$, a flow of $|\mathcal{A}_b|$ leaves entry point $b$, and flow is conserved in all other nodes of the graph.
Constraint~\eqref{eq:flow-u} enforces edges with a positive flow to participate in the arrival route, where $Q$ is a large number (e.g., $Q=|\mathcal{A}|$).
Constraint~\eqref{eq:mcf-indegree-vstar} ensures one ingoing edge for the runway $r$, and  Constraint~\eqref{eq:mcf-outdegree-S} guarantees that each entry point has one outgoing edge. 
For all other nodes, Constraints~\eqref{eq:mcf-indegree} and~\eqref{eq:mcf-outdegree} ensure  the maximum indegree of two and outdegree of one, respectively.
Constraint~\eqref{eq:curvature} prevents sharp turns by enforcing that either  edge $(i,j)$ is used (in which case no edge in 
$\Gamma_{i,j}$ may be used), or any subset of edges in $\Gamma_{i,j}$ may be chosen.
Although degree constraints prevent routes from crossing at nodes, routes crossing may occur within a grid square when additional constraints such as temporal separations are included.
Routes crossings within a grid square are prevented by Constraints~\eqref{eq:aux-1}--\eqref{eq:aux-5}, where $V' = V\setminus\{\mbox{last grid row}\}\setminus\{\mbox{last grid column}\}$.

Variables ${x}^{\scalebox{0.7}{EP}}_{(i,j),b}$ and $x_{i,j}$ are connected by Constraint~\eqref{eq:xbs1}.
Constraints~\eqref{eq:xbs2} and~\eqref{eq:xbs3} ensure that exactly one path from each entry point goes to the runway and Constraint~\eqref{eq:xbs4} enforces that the flow   is conserved in all the other nodes.

Constraint~\eqref{eq:y1n_2} ensures that for each aircraft $a$ from its entry point $b$, exactly one speed profile $p$ is chosen and the aircraft arrives to its entry point at exactly one time in $[\bar{t}_a-\mu, \bar{t}_a+\mu]$. 
Several of the $y$-variables  are set to zero in Constraints~\eqref{eq:y2n_2}--\eqref{eq:nowhere-before}, that is, whenever we know that an aircraft cannot occupy the node at all or not at certain points in time. 
Constraint~\eqref{eq:y6n} combines the $y$-variables and the $x_{e,b}$ and ensures that any aircraft $a$ can occupy a node $j$ at any time $t$ only if there exists an ingoing edge for $j$.
Each aircraft uses maximum one speed profile at any time and node by Constraint~\eqref{eq:y8n}. 

For each aircraft $a$ arriving at entry point $b$, the selected speed profile from $\mathcal{S}(a)$  must have a length equal to  $\ell(b)$, which is determined by the arrival path from $b$.
In particular, Constraint~\eqref{eq:length} calculates $\ell(b)$ by counting the $x^{\scalebox{0.7}{EP}}_{(i,j),b}$ variables, while Constraint~\eqref{eq:length-profile} sets the binary variable  $y$ to 1 for the speed profile whose length matches $\ell(b)$. 
Finally, Constraints~\eqref{eq:psi1}--\eqref{eq:psi4} enforce that only the speed profile with the correct length  $\ell(b)$ is selected for each aircraft.

For all categories $C_{\kappa_1}$ and $C_{\kappa_2}$, a temporal separation of $\sigma_{\kappa_1, \kappa_2}$ is enforced by using Constraint~\eqref{eq:wake-turb-1}, and if $\sigma_{\kappa_1, \kappa_2}\geq 2$ then Constraint~\eqref{eq:wake-turb-1'} is also used.
Similarly, for all categories $C_{\kappa}$,  a temporal separation of $\sigma_{\kappa, \kappa}$  is applied by Constraint~\eqref{eq:wake-turb-2}, and if $\sigma_{\kappa, \kappa}\geq 2$ then Constraint~\eqref{eq:wake-turb-2'} is applied as well.

An aircraft $a$ along its route using speed profile $p$ arrives in the $k+1$-st node $i$ at time $t$ only if it traverses an edge $(j,i)$ which takes $u_{a,p,k}$, that is, the aircraft is on node $j$ at time $t-u_{a,p,k}$. 
These times at  which an aircraft arrives at  nodes along its arrival route are given by Constraints~\eqref{eq:z1n}--\eqref{eq:z3n}. Finally, the domains of variables are given in~\eqref{eq:dom1}--\eqref{eq:dom8}

\subsection{Path-Based Formulation}\label{subsec:new-model}
To reformulate the problem as a path-based model, we first describe the components of the model included in the definition of the paths and a method for constructing these paths. 
In Subsection~\ref{subsec:gen-paths}, we detail the process of generating all possible 
paths from each entry point to the runway, along with the trajectories of each aircraft on these paths based on their CDO speed profiles. 
Thereafter, we present the reformulated model that utilizes these pre-generated paths in Subsection~\ref{subsec:reform}, 
followed by a discussion of some special cases of the problem in Subsection~\ref{subsec:special}.

\subsubsection{Generating Paths and Trajectories}\label{subsec:gen-paths}
A key step in formulating a path-based model is defining the set of feasible paths from each entry point to the runway while ensuring that they follow operational constraints.
These paths must be generated with respect to forbidden sharp turns and obstacle avoidance.
Obstacle avoidance can be handled easily by removing  the nodes and edges of that area from the graph. 

\begin{figure*}[h!] \hspace*{-1cm}
	\begin{subfigure}[t]{0.4\textwidth}
	\centering
	\includegraphics[width=0.5\textwidth]{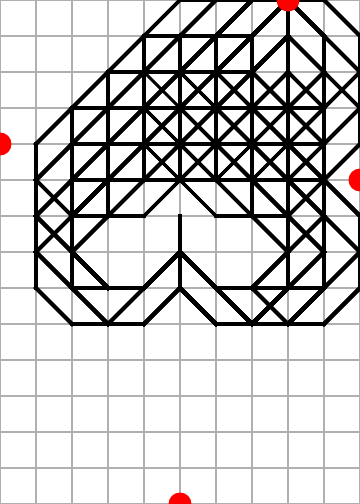} \caption{HMR}
\end{subfigure}\hfill \hspace*{-2.75cm}
\begin{subfigure}[t]{0.4\textwidth}
	\centering
	\includegraphics[width=.5\textwidth]{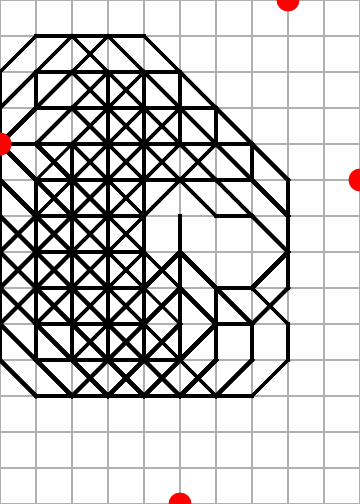} \caption{ELTOK}
\end{subfigure}\hfill \hspace*{-2.75cm}
\begin{subfigure}[t]{0.4\textwidth}
	\centering
	\includegraphics[width=.5\textwidth]{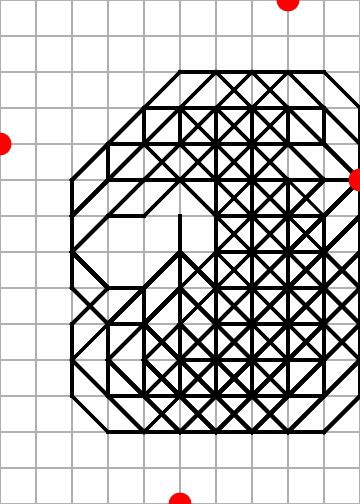} \caption{XILAN}
\end{subfigure}\hfill \hspace*{-2.75cm}
\begin{subfigure}[t]{0.4\textwidth}
	\centering
	\includegraphics[width=.5\textwidth]{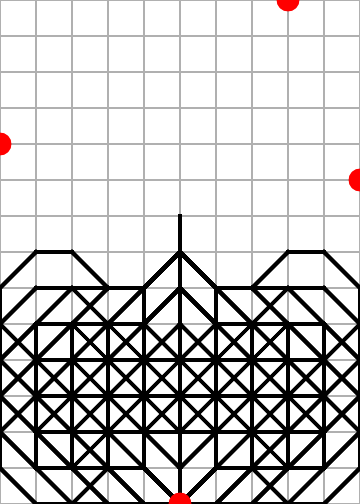} \caption{NILUG}
\end{subfigure}\hspace*{-1cm}
\caption{The areas covered by length-bounded paths with upper bound of 14 from each entry point of Arlanda TMA}\label{fig:arlanda-TMA-path14} 
\end{figure*}

\begin{figure*}[h!] \hspace*{-1cm}
	\begin{subfigure}[t]{0.4\textwidth}
	\centering
	\includegraphics[width=0.55\textwidth]{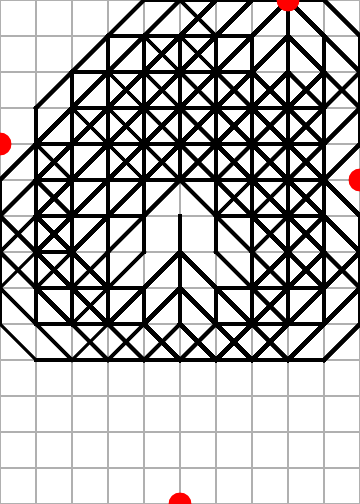} \caption{HMR}
\end{subfigure}\hfill \hspace*{-2.75cm}
\begin{subfigure}[t]{0.4\textwidth}
	\centering
	\includegraphics[width=.5\textwidth]{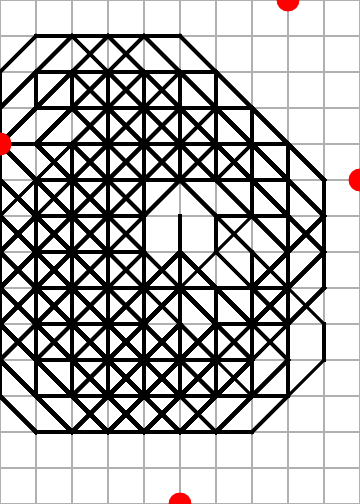} \caption{ELTOK}
\end{subfigure}\hfill \hspace*{-2.75cm}
\begin{subfigure}[t]{0.4\textwidth}
	\centering
	\includegraphics[width=.5\textwidth]{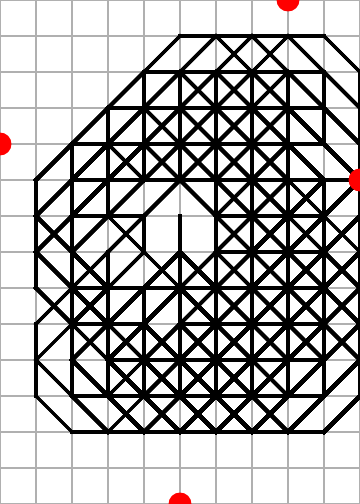} \caption{XILAN}
\end{subfigure}\hfill \hspace*{-2.75cm}
\begin{subfigure}[t]{0.4\textwidth}
	\centering
	\includegraphics[width=.5\textwidth]{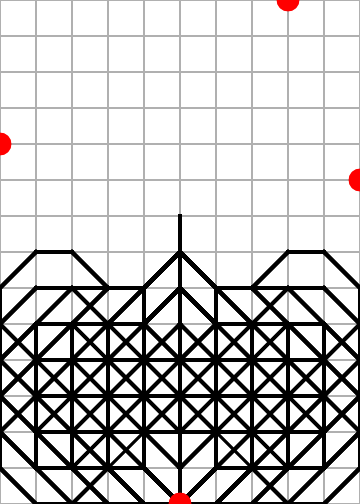} \caption{NILUG}
\end{subfigure}\hspace*{-1cm}
\caption{The areas covered by length-bounded paths with upper bound of 15 from each entry point of Arlanda TMA}\label{fig:arlanda-TMA-path15} 
\end{figure*}

Additionally, we consider an upper bound on the length of the paths. 
Similar to Subsection~\ref{subsec:grid-turbulence}, let $\lambda$ denote the upper bound on the number of nodes in any path, that is, any path will include maximum $\lambda-1$ edges.
An upper bound on path length is necessary to balance computational efficiency and practical feasibility. 
From each entry point, the nodes reachable within this limit cover a large portion of the TMA, ensuring sufficient routing options. 
For instance, the areas covered by paths from each entry point in the TMA of Stockholm Arlanda airport with upper bounds of 14 and 15 are given in Figure~\ref{fig:arlanda-TMA-path14}  and~\ref{fig:arlanda-TMA-path15}, respectively.
Moreover, since the objective function is minimizing both tree size and paths length, excessively long paths would not be chosen.
Therefore, generating them would be redundant and would unnecessarily complicate the model.

In order to present the method, we introduce additional notation. 
Assume that for each node $i$, $\Lambda_i$ is the set of successors of node $i$, that is, $\{j\in V: (i,j)\in E\}$.
If edge $e=(i,j)$ is used in the arrival tree, all the outgoing edges from $j$ must form an angle of at least $\gamma$ with $e$. 
Let $\bar{\Gamma}_e$ be the set of all nodes $\bar{j}\in \Lambda_j$  such that  edge $\bar{e}=(j, \bar{j})$ 
forms an angle smaller than $\gamma$ with $e$. 
Therefore, if $e$ is used in the arrival tree, none of the nodes from $\bar{\Gamma}_e$ will be visited after $j$ in the tree, 
that is, edge  $\bar{e}$ cannot be part of the tree. 
Moreover, for implementation purposes, let $\mathbf{\Lambda}$ be the adjacency list of the graph 
and $\mathbf{\bar{\Gamma}}$ be a list  of all $\bar{\Gamma}_e$.

To generate length-bounded paths, we first determine the minimum number of edges needed to reach the runway from each node. 
This information helps us enforce the path-length bound efficiently. 
From each node $i\in V$, we compute these minimum values, denoted by $\delta_i$, 
using a shortest path algorithm based on  breadth-first search (BFS)  applied in reverse order, starting from the runway.
Let $\mathbf{\delta}=(\delta_i)_{i\in V}$ be these values for all nodes in the graph.

Once $\mathbf{\delta}$ is available, all possible paths from a start node $\bar{i}$ to an end node, the runway $r$ in our case, 
are calculated by Algorithm \ref{alg-1}, where $\iota$, initialized as zero,  is a counter for controling the number of edges in the paths and
\textbf{find\_all\_paths} is a recursive function given in Algorithm \ref{alg-2}.
This function, detailed below, essentially performs a depth-first search (DFS) to find all paths from the start node to the end node, 
while ensuring that the paths comply with the upper bound on path length and forbidden sharp turns.

\vspace*{0.5cm}
\begin{algorithm}[b]
	\caption{Find all paths from node $\bar{i}$ to $r$}\label{alg-1}
	\begin{algorithmic}
		\State \textbf{Input:} $\mathbf{\Lambda}$, $\bar{i}$, $r$, $\lambda$, $\iota=0$, $\mathbf{\bar{\Gamma}}$, $\mathbf{\delta}$, $\pi=[~]$.
		\State \textbf{Output:} $\Pi_{\bar{i}}=$ \textbf{find\_all\_paths}($\bar{i}$, $\iota=0$, $\pi=[~]$)
	\end{algorithmic}
\end{algorithm}
\vspace*{0.5cm}
\begin{algorithm}[b]
	\caption{\textbf{find\_all\_paths}($\bar{i}$, $\iota$, $\pi$)}\label{alg-2}
	\begin{algorithmic}
		\State Add node $\bar{i}$ to the end of path $\pi$.
		\State Let ${\bar{e}}$ be the last edge of $\pi$ (If $\bar{i}$ is the first node of $\pi$, let $\bar{e}=0$ and  $\bar{\Gamma}_{0}=\emptyset$).
		\If{$\bar{i}=r$}
		\State \textbf{return} $\Pi_{\bar{i}}=\{\pi\}$
		\EndIf
		\State  $\Pi_{\bar{i}} = [~]$
			\For{ $\hat{i}\in\Lambda_{\bar{i}}$}
			\If{ $\hat{i}\notin\pi$ and $(\iota + \delta_{\hat{i}}) < (\lambda - 1)$ and $\hat{i}\notin \bar{\Gamma}_{{\bar{e}}}$}
			\State $\iota = \iota + 1$
			\State   $\hat{\Pi}$ = \textbf{find\_all\_paths}($\hat{i}$, $\iota$, $\pi$) 
			\For{$\hat{\pi} \in \hat{\Pi}$}
			\State Add   $\hat{\pi}$ to    $\Pi_{\bar{i}}$
			\EndFor
			\State $\iota = \iota - 1$
			\EndIf
			
			\EndFor
		\State \textbf{return} $\Pi_{\bar{i}}$
	\end{algorithmic}
\end{algorithm}

The function starts from an empty path $\pi$, adds the start node to it and determines the last edge of the current path. 
If the current node is the end node, it returns the current path as the only path in $\Pi_{\bar{i}}$ and stops. 
Otherwise, it considers the neighboring nodes of the current node, that is, $\Lambda_{\bar{i}}$,  and for each neighboring node $\hat{i}$, it checks three conditions: 
It is not already in the current path, the current path length plus the shortest distance to $r$ is less than the upper bound $ub$, and it does not make a sharp turn. 
If all conditions are satisfied, it increments the path length counter and recursively calls the function to find all paths from $\hat{i}$ to $r$. 
Then, each path found from the recursive call is added to the list of all valid paths $\Pi_{\bar{i}}$ and the path length counter is decremented after processing.
In the end, the list of all valid paths $\Pi_{\bar{i}}$ is returned. 

For each entry point, we obtain all feasible paths using Algorithm \ref{alg-1}. 
Let $\Pi_b$ denote the set of paths from entry point $b$ to the runway and 
$\mathbf{\Pi}=\bigcup_{b\in \mathcal{P}}\Pi_b$ be the set of all paths. 
Each path is an ordered list of nodes and one can extract different information from this list and/or generate the required data based on the paths. 
In the rest of this subsection,  we present the data extracted from the paths to be used 
in the model in Subsection \ref{subsec:reform}.

For any $\pi\in\mathbf{\Pi}$, we define $\theta_{\pi}$ as the set of edges that path $\pi$ passes through. 
Moreover,  we denote the set of all paths crossing node $i\in V$ by $\overline{\Pi}_i$ and the set of all paths going through edge $(i,j)\in E$ by $\zeta_{i,j}$.
Since the length of each path is known, we can specify the optimal speed profile of each aircraft on each path. 
To define the trajectories, we first note that each aircraft can only enter at its entry point and can therefore only use paths defined for this entry point. 
In the case of a fixed entry time, there is a one-to-one correspondence between a trajectory of an aircraft and a path from that entry point. 
In the case of flexible entry times, one trajectory is created for each of the possible entry times in the interval  $[\bar{t}_a-\mu,\bar{t}_a+\mu]$, that is, 
$2\mu + 1$ trajectories for each aircraft, assuming a discretization of entry times to minutes.

The time when aircraft $a$ with entry time $\bar{t}_a$ occupies node $i$ on path $\pi$ based on pre-computed CDO speed profiles is denoted by
$\xi_{a,\pi,i}$. We define two sets based on $\xi_{a,\pi,i}$. 
For an aircraft $a$, node $i$ and time $t$, $ \Xi_{a,i,t}$ is defined as 
the set of all  paths  on which aircraft $a$ with entry time $\bar{t}_a$ occupies node $i$ at time $t$, that is, 
$\Xi_{a,i,t} = \{\pi\in\overline{\Pi}_i\cap\Pi_b: ~~ \xi_{a,\pi,i} = t\}.$
Moreover, for an aircraft $a$, node $i$, time $t$ and temporal separation $\sigma$, we define  $\Upsilon_{a,i,t,\sigma}$ as 
the set of all  paths  on which aircraft $a$ with entry time $\bar{t}_a$ occupies node $i$ sometime between time $t$ and  $t+\sigma-1$, that is, 
$\Upsilon_{a,i,t,\sigma} = \{\pi\in\overline{\Pi}_i\cap\Pi_b: ~~ t\leq \xi_{a,\pi,i}\leq t+\sigma-1\}.$

\subsubsection{Reformulated Model}\label{subsec:reform}
Given all the sets and parameters defined in Subsections~\ref{subsec:grid-turbulence} and \ref{subsec:gen-paths}, 
we present the path-based formulation of the complete model given by~\cite{spshspp-asmda-21}.
We introduce binary variables
\begin{itemize}
	\item  $\rho_{\pi}$  indicates whether or not path $\pi$ is used in the arrival tree.
	\item  $x_{i,j}$ indicates whether  or not edge $(i,j)\in E$ participates in the arrival tree. 
	\item  To handle flexible entry times, $\tau_{a,\pi,t}$ gives
	whether or not aircraft $a$ uses path $\pi$ and arrives at its corresponding entry point at time $t$. 
\end{itemize}
We call the resulting model M2:




\setcounter{equation}{-1}
\makeatletter
\renewcommand{\theequation}{M2.\arabic{equation}}
\makeatother

 {\allowdisplaybreaks
 	\footnotesize
\begin{align}\label{obj}
	\min~~~ &\beta\sum_{(i,j)\in E} l_{i,j}x_{i,j} +(1-\beta)\sum_{b \in \mathcal{P}}\sum_{\pi\in\Pi_b}\sum_{(i,j)\in \theta_{\pi}} |\mathcal{A}_b|l_{i,j}\rho_{\pi}\\[.5em]
	\label{indeg}
	\text{s.t.}~~~&\sum_{j:(j,i)\in E} x_{j,i} \leq 2,&&\hspace*{-1cm}\forall i\in V\setminus\{\mathcal{P}\cup r\},\\\label{outdeg}
	&\sum_{j:(i,j)\in E}x_{i,j}\leq 1,&&\hspace*{-1cm}\forall i\in V\setminus\{\mathcal{P}\cup r\},\\
	&x_{i,i+1+n}+x_{i+1+n,i}+x_{i+n,i+1}+x_{i+1,i+n}\leq 1,&&\nonumber\\
	&&&\hspace*{-2.5cm}\forall i\in V'\setminus\{\mathcal{P}\cup r\}: 
i+1+n, i+n, i+1\notin \{\mathcal{P}\cup r\},\label{cross1}\\
	&x_{i,i+1+n}+x_{i+n,i+1}+x_{i+1,i+n}\leq 1,  &&\hspace*{-1cm}\forall i \in \mathcal{P}\cap V',\label{cross2}\\
	&x_{i,i+1+n}+x_{i+1+n,i}+x_{i+1,i+n}\leq 1,  &&\hspace*{-1cm}\forall i: i+1 \in \mathcal{P},\label{cross3}\\
	&x_{i,i+1+n,}+x_{i+n+1,i}+x_{i+n,i+1}\leq 1, &&\hspace*{-1cm}\forall i: i+n \in \mathcal{P},\label{cross4}\\
	&x_{i+1+n,i}+x_{i+n,i+1}+x_{i+1,i+n}~\leq 1,  &&\hspace*{-1cm}\forall i: i+n+1 \in \mathcal{P}, \label{cross5}\\
	&\sum_{\pi\in\Pi_b}\rho_{\pi}=1,&&\hspace*{-1cm}\forall b\in\mathcal{P},\label{onepath}\\
	\label{coup2}
	&\rho_{\pi}\leq x_{i,j},&&\hspace*{-1cm}\forall\pi\in\mathbf{\Pi}, \forall (i,j)\in \theta_{\pi},\\\label{tra1}
	&\sum_{\pi \in \Pi_b}\sum_{t=\bar{t}_a-\mu}^{\bar{t}_a+\mu} \tau_{a,\pi,t} = 1,&&\hspace*{-1cm}\forall b\in \mathcal{P}, \forall a \in \mathcal{A}_b,\\\label{tra2}
	&\sum_{t=\bar{t}_a-\mu}^{\bar{t}_a+\mu} \tau_{a,\pi,t} = \rho_{\pi},&&\hspace*{-1cm}\forall b\in \mathcal{P}, \forall a \in \mathcal{A}_b, \forall \pi \in \Pi_b,\\
	&\sum_{b \in\mathcal{P}}\sum_{a \in \mathcal{A}_b\cap C_{\kappa_2}}\sum_{t'=-\mu}^{\mu} \hspace*{-0.6cm}\sum_{\hspace*{.7cm}\pi\in\Upsilon_{a,i,t-t',\sigma_{\kappa_1,\kappa_2}}} \hspace*{-0.5cm}\tau_{a,\pi,\bar{t}_a+t'}\leq&&\nonumber\\
	&&&\hspace*{-2.5cm}\Omega(1-\sum_{b \in\mathcal{P}}\sum_{a' \in \mathcal{A}_b\cap C_{\kappa_1}}\sum_{t'=-\mu}^{\mu} \hspace*{-0.6cm}\sum_{\hspace*{.7cm}\pi'\in\Xi_{a',i,t-t'}}\hspace*{-0.5cm} \tau_{a',\pi',\bar{t}_{a'}+t'}),\nonumber\\
	&&&\hspace*{-3.6cm}\forall \kappa_1, \kappa_2\in \mathcal{C}: \kappa_1\neq \kappa_2, ~\forall i\in V, ~\forall t\in\{0,\ldots,\overline{T}-\sigma_{\kappa_1, \kappa_2}+1\},	\label{temp-diff}\\
	\nonumber
	&\sum_{b \in\mathcal{P}}\sum_{a \in \mathcal{A}_b\cap C_{\kappa}:  a\neq a' }\sum_{t'=-\mu}^{\mu} \hspace*{-0.6cm}\sum_{\hspace*{.7cm}\pi\in\Upsilon_{a,i,t-t',\sigma_{\kappa,\kappa}}} \hspace*{-0.5cm}\tau_{a,\pi,\bar{t}_a+t'}\leq &&\hspace*{-0.4cm}\Omega(1-\sum_{t'=-\mu}^{\mu} \hspace*{-0.6cm}\sum_{\hspace*{.7cm}\pi'\in\Xi_{a',i,t-t'}} \hspace*{-0.5cm}\tau_{a',\pi',\bar{t}_{a'}+t'}),\hspace*{3cm}\\\label{temp-sam}
	&&&\hspace*{-3.6cm}\forall \kappa\in\mathcal{C}, \forall a'\in C_{\kappa}, \forall i\in V, ~
\forall t\in\{0,\ldots,\overline{T}-\sigma_{\kappa, \kappa}+1\},\\\label{dom-1}
	&x_{i,j}~~~~\text{binary},  &&\hspace*{-1cm}\forall (i,j)\in E,\\\label{dom-2}
	&\rho_{\pi}~~~~~~\text{binary},&&\hspace*{-1cm}\forall\pi\in\mathbf{\Pi},\\\label{dom-3}
	&\tau_{a,\pi,t}~~\text{binary},&&\hspace*{-1cm}\forall b\in \mathcal{P}, \forall a \in \mathcal{A}_b, 
\forall t \in [\bar{t}_a-\mu, \bar{t}_a+\mu].
\end{align}}

The objective function~\eqref{obj} remains structurally the same as in the original model~\eqref{obj-1}, consisting of a convex combination of the weight of the arrival-route tree  and the total aircraft distances along the paths.
However, in the reformulated model, the second term (the total aircraft distances along the paths) is expressed in terms of the paths utilized. 
Specifically, the selected path for entry point $b$ is traversed $|\mathcal{A}_b|$ times.
Therefore, the total distances of all aircraft originating from a specific entry point is given by the number of aircraft arriving at that entry point multiplied by the length of the selected path from that entry point (which effectively prioritizes paths serving higher demand).

Constraints~\eqref{indeg} and~\eqref{outdeg} correspond directly to~\eqref{eq:mcf-indegree} and~\eqref{eq:mcf-outdegree}, ensuring the required indegree and outdegree conditions are satisfied.
Preventing crossings is achieved by Constraints~\eqref{cross1}--\eqref{cross5} which mirror Constraints~\eqref{eq:aux-1}--\eqref{eq:aux-5} in the original model.
Alternatively, for each path, we can identify the set of paths that cross it during path generation, 
and then add a constraint ensuring that the current path and any path from its crossing set cannot be selected at the same time.  
However, our initial tests showed that this approach resulted in worse computational performance, so we did not include these constraints in our final model.

Constraint~\eqref{onepath} states that exactly one path from each entry point must be used. 
This constraint replaces Constraints~\eqref{eq:mcf-indegree-vstar}, \eqref{eq:mcf-outdegree-S}, \eqref{eq:xbs2} and~\eqref{eq:xbs3} as each path originates from an entry point and terminates at the runway.
Constraint~\eqref{coup2} establishes connections between variables, effectively replacing Constraint~\eqref{eq:xbs1}.  
Constraint~\eqref{tra1} ensures that for each aircraft, exactly one trajectory is selected and represents a reformulation of Constraint~\eqref{eq:y1n_2}.
Constraint~\eqref{tra2} states that if a path is not used, none of the trajectories on that path can be used and if a path is used, 
one trajectory for each aircraft entering the corresponding entry point on that path must be used. 
Note that Constraint~\eqref{tra2} together with~\eqref{onepath} yield Constraint~\eqref{tra1}. The initial tests showed that adding these constraints improves the computational performance. Therefore, we include them in the model. 

Temporal separation for different wake-turbulence categories is ensured by Constraints~\eqref{temp-diff}--\eqref{temp-sam}:
Constraint~\eqref{temp-diff} enforces the temporal separation when the two aircraft are from different categories, 
and  Constraint~\eqref{temp-sam} enforces the temporal separation  when both aircraft belong to the same category. 
These constraints are reformulations of Constraints~\eqref{eq:wake-turb-1}--\eqref{eq:wake-turb-2'}, expressed in terms of the new variables. 
Finally,  Constraints~\eqref{dom-1}--\eqref{dom-3} define the domains of variables.

The  model presented in Subsection~\ref{subsec:old-model} is based on arcs and flows on arcs while our reformulated model is based on paths. 
Variables $x_{i,j}$ and Constraints~\eqref{indeg}, \eqref{outdeg}, and~\eqref{cross1}--\eqref{cross5} remain unchanged from the compact model, ensuring consistency in structure where applicable. 
Other constraints, however, are either reformulated using the new variables or incorporated during the path and trajectory generation process.
In the reformulated model, paths are carefully generated to avoid sharp turns, ensuring that only feasible routes are considered. 
For each generated path, the length is fixed, allowing us to associate the appropriate speed profile with the path for various entry times. 
Consequently, the trajectories are predefined, and the remaining task is to select the appropriate trajectory, which is managed by the variables $\tau_{a,\pi,t}$. 
As a result, all constraints in the compact model related to speed profiles are effectively handled during the trajectory generation, as described in Subsection~\ref{subsec:gen-paths}. 
Finally, this reformulation has generated all feasible binary points of the Dantzig-Wolfe discretization with respect to constraints~\eqref{eq:curvature}, \eqref{eq:xbs4}, \eqref{eq:y2n_2}--\eqref{eq:psi4}, and~\eqref{eq:z1n}--\eqref{eq:yaktu-new}, represented by the paths and trajectories in our approach.

\subsubsection{Special Cases}\label{subsec:special}
A simplified version of the problem arises when entry times are fixed and/or only a single aircraft category is considered. 
With fixed entry times, the scheduling becomes less complex, as the timing of each aircraft's arrival is predetermined. 
However, this may lead to infeasible solutions if the fixed times conflict with temporal separation requirements.
In the case of a single aircraft category, the temporal separation requirement, denoted by 
$\sigma$, is uniform across all aircraft, simplifying the problem by eliminating the need for different separation criteria based on aircraft type.
In the following, we discuss how our model is modified to address the simplified version of the problem.

If flexible entry times are allowed, but there is only one aircraft category,
Constraints~\eqref{temp-diff}--\eqref{temp-sam} are replaced by Constraint~\eqref{tra3}
{\small
\begin{align}
	&\sum_{b \in\mathcal{P}} \sum_{a\in\mathcal{A}_b}\sum_{t'=-\mu}^{\mu} \hspace*{-0.6cm}\sum_{\hspace*{.6cm}\pi\in\Upsilon_{a,i,t-t',\sigma}} \hspace*{-0.3cm}\tau_{a,\pi, \bar{t}_a+t'}\leq 1,
	 &&\forall i\in V,~\forall t\in\{0,\ldots,\overline{T}-\sigma\}.\label{tra3}	
\end{align}}
If the entry times are fixed and there is only one aircraft category, variables $\tau_{a\pi t}$ and Constraints~\eqref{tra1}--\eqref{temp-sam} are not needed.
Instead, we add Constraint~\eqref{tempo} to ensure the temporal separation
{\small
\begin{align}
		&\sum_{b \in\mathcal{P}} \sum_{a\in\mathcal{A}_b} \sum_{\pi\in\Upsilon_{a,i,t,\sigma}} \rho_{\pi}\leq 1,
		&&\forall i\in V,~\forall t\in\{0,\ldots,\overline{T}-\sigma\}.\label{tempo}	
\end{align}	}

\subsection{Tree Consistency}\label{subsec:consistency}
The arrival trees are optimized based on the aircraft currently arriving, as traffic at different times may require routes of varying lengths. 
Therefore, these trees need to be recomputed frequently to reflect the traffic changes in the TMA. 
However, frequent updates can lead to instability, as large differences between consecutive trees may complicate the management of air traffic.
To address this, we aim to maintain consistency by limiting the number of edges that differ between trees for consecutive time periods.

We measure this consistency using edge indicator variables $x_{i,j}$ for the current time period and $x_{i,j}^{old}$ for the previous one. 
Let $E_1$ denotes the set of edges included in the previous tree, defined as $E_1=\{(i,j)\in E: x_{i,j}^{old}=1\}$, and let $E_0=E\setminus E_1$ represents the remaining edges.
A new variable, $x_{i,j}^{\Delta}$, is introduced to denote $|x_{i,j}-x_{i,j}^{old}|$, as specified in Constraints~\eqref{eq:consist-1}--\eqref{eq:consist-2}.
Finally, to limit the variation between consecutive trees, we enforce Constraint~\eqref{eq:consist-3}, which bounds the number of differing edges using the parameter $U$.

{\small
\begin{align}
	&x_{i,j}^{\Delta} = 1-x_{i,j},&& \forall (i,j)\in E_1\label{eq:consist-1}\\
	& x_{i,j}^{\Delta} = x_{i,j}, && \forall (i,j)\in E_0\label{eq:consist-2}\\
	&\sum_{(i,j)\in E} x_{i,j}^{\Delta} \leq U, &&\label{eq:consist-3}
\end{align}}

To fulfill Constraint~\eqref{eq:consist-3}, the model may retain certain edges from the previous tree even if they are not part of any current path. 
To prevent this behavior, we introduce the following constraint, which ensures that an edge is included in the solution only if it is traversed by at least one selected path.

{\small
\begin{align}
	&x_{i,j}\leq \sum_{\pi: (i,j)\in\theta_{\pi}}\rho_{\pi},&& \forall (i,j)\in E_1.\label{eq:consist-4}
\end{align}}

\subsection{Presence of Aircraft from the Previous Period in the TMA}\label{subsec:aircraft-TMA}
Assume that a tree from the previous time period is available, and that some aircraft have not yet reached the runway by the beginning of the current period.
Consequently, at the start of the new period, certain aircraft remain in TMA.
For safety reasons, these aircraft must be taken into account when designing the tree for the current period.
To achieve this, we fix the positions of these aircraft and enforce temporal separation from them.
Suppose that aircraft $a'\in C_{\kappa_1}$ is located at node $i$ at time $t$ where $t$ is within the current time period (i.e., after the start of the period).
Aircraft that have already reached the runway by this time are excluded from consideration.
To ensure proper temporal separation, no other aircraft from $C_{\kappa_2}$ may occupy node $i$ within time $t$ and $t+\sigma_{\kappa_1, \kappa_2}-1$
and within time $t-\sigma_{\kappa_1, \kappa_2}+1$ and $t$, for any $\kappa_2\in \mathcal{C}$. 
Accordingly, for each aircraft $a'\in C_{\kappa_1}$ located at node $i$ at time $t$, we enforce Constraints~\eqref{eq:airATM-1} and~\eqref{eq:airATM-2}.
{\small
\begin{align}
	&\sum_{b \in\mathcal{P}}\sum_{a \in \mathcal{A}_b\cap C_{\kappa_2}}\sum_{t'=-\mu}^{\mu} \hspace*{-0.6cm}\sum_{\hspace*{.7cm}\pi\in\Upsilon_{a,i,t-t',\sigma_{\kappa_1,\kappa_2}}} \hspace*{-0.5cm}\tau_{a,\pi,\bar{t}_a+t'}=0,&&\forall  \kappa_2\in \mathcal{C},\label{eq:airATM-1}\\
&\sum_{b \in\mathcal{P}}\sum_{a \in \mathcal{A}_b\cap C_{\kappa_2}}\sum_{t'=-\mu}^{\mu} \hspace*{-0.6cm}\sum_{\hspace*{.7cm}\pi\in\Upsilon_{a,i,t-t'-\sigma_{\kappa_2,\kappa_1}+1,\sigma_{\kappa_2,\kappa_1}}} \hspace*{-1.2cm}\tau_{a,\pi,\bar{t}_a+t'}=0,&&\forall  \kappa_2\in \mathcal{C}.\label{eq:airATM-2}
\end{align}}

Moreover, if an aircraft is traveling along edge  $(i,j)$ from time $t_1$ to time $t_2$, then no other aircraft is allowed to occupy the reverse edge $(j,i)$ or any edge that crosses $(i,j)$ during the same time interval.
Suppose aircraft $a'$ is on edge $(i,j)$ within time $t_1$ to time $t_2$.
If $(i,j)$ is a diagonal edge of the grid, then  $(i', j')$ and $(j', i')$ denote  the crossing edges. 
In this case, all these conflicting edges must remain unused by other aircraft during the interval. To ensure this, we consider Constraints~\eqref{eq:airATM-3}--\eqref{eq:airATM-5}.

{\small
\begin{align}
	&\sum_{b \in\mathcal{P}}\sum_{a \in \mathcal{A}_b}\sum_{t'=-\mu}^{\mu} \sum_{t=t_1}^{t_2}\hspace*{-0.6cm}\sum_{\hspace*{.6cm}\pi\in\Xi_{a,j,t}\cap \zeta_{j,i}} \hspace*{-0.6cm}\tau_{a,\pi,\bar{t}_a+t'}=0,&&\label{eq:airATM-3}\\
	&\sum_{b \in\mathcal{P}}\sum_{a \in \mathcal{A}_b}\sum_{t'=-\mu}^{\mu} \sum_{t=t_1}^{t_2}\hspace*{-0.6cm}\sum_{\hspace*{.6cm}\pi\in\Xi_{a,j',t}\cap \zeta_{j',i'}} \hspace*{-0.6cm}\tau_{a,\pi,\bar{t}_a+t'}=0,&&\label{eq:airATM-4}\\
	&\sum_{b \in\mathcal{P}}\sum_{a \in \mathcal{A}_b}\sum_{t'=-\mu}^{\mu} \sum_{t=t_1}^{t_2}\hspace*{-0.6cm}\sum_{\hspace*{.6cm}\pi\in\Xi_{a,i',t}\cap \zeta_{i',j'}} \hspace*{-0.6cm}\tau_{a,\pi,\bar{t}_a+t'}=0,&&\label{eq:airATM-5}
\end{align}}

%% file: comput-res.tex
\section{Computational Results}\label{comp-sec}
In this section, we present the computational results of our extended model. 
The implementation of generating paths and trajectories was done in Python, 
and we solved our model using Gurobi version 9.5.1 as the solver with gurobipy as the interface to Python version 3.10.4.
We ran the experiments on a MacBook Pro, M1 2020.

We used real-world data of arriving traffic at Stockholm Arlanda TMA with a data sample for one of the busiest hours of
 operation in 2018: May 16, from 5:00 AM to 5:59 AM, with 33 aircraft, and the following hour, from 6:00 AM to 6:59 AM, with 16 aircraft.
The data was obtained from the open-source database of the Opensky Network (\cite{opensky}).
The data set contains aircraft arriving to runway 19L (one of the three runways of Arlanda airport, where 01L/19R and 01R/19L are parallel), and the four main entry points HMR (north), XILAN (east), NILUG (south), and ELTOK (west).
Moreover, we use an $11\times 15$ grid  which ensures a separation of 6 NM.
Arlanda TMA with the runway, entry points, and the overlaid grid  is illustrated in Figure~\ref{fig:arlanda-TMA-grid}(b), which is given by~\cite{spshspp-asmda-21}.

First, we explain the parameters' setup, 
details of our model, and the path generation algorithm in Subsection~\ref{subsec:param}.  
Then, we discuss two experiment sets in Subsection~\ref{subsec:exper}.
Each experiment is conducted with different scenarios and 
in case of only one aircraft category,  we replace Constraints~\eqref{temp-diff}--\eqref{temp-sam} with Constraint~\eqref{tra3}.

\subsection{Parameters' Setup}\label{subsec:param}
We generated all possible paths from each entry point to the runway using Algorithm~\ref{alg-1} for two bounds on path length: 14 and 15 edges.
The  minimum angle of the outgoing edge was set to 135 degrees.
The running time to generate all paths with a maximum of 14 edges was 0.22 seconds, while for 15 edges, it was 0.52 seconds.
We obtained a total of 3821 and 9527 paths from all entry points to the runway for an upper bound of 14 and 15 edges, respectively.

In order to generate the trajectories on the paths, we need speed profiles. 
The CDO speed profiles determine the time each aircraft spends to traverse each TMA segment. 
For each aircraft,  speed profiles for different route lengths were computed by~\cite{spshspp-asmda-21}. 
We use the same speed profiles in our experiments.

Moreover, we use aircraft categories from \cite{icao-cdom-10}: 
LIGHT (L), MEDIUM (M), HEAVY (H). 
We deﬁne $C_1 = \{H,M\}$, and $C_2 = \{L\}$.
Based on ICAO's separation minima, see \cite{icao-wt-20}, we let
$\sigma_{1,2}=3$ minutes, and  $\sigma_{1,1}=\sigma_{2,2}=\sigma_{2,1}=2$ minutes, that is,
3 minutes of temporal separation if the trailing aircraft is light and the leading aircraft is medium or heavy,  and 2 minutes otherwise.
In the case of only one category, the temporal separation is 2 minutes,  $\sigma=2$.
Finally, we let  $\beta=0.1$ in the objective function, that is, we prioritize the weighted shortest paths in the objective function.

\subsection{Experiments}\label{subsec:exper}
We run two sets of experiments and  in each experiment, we consider different scenarios where some additional constraints such as tree consistency and/or presence of aircraft from the privous period in the TMA are included. 
In each case, we run it once with the aircraft present in the real-world data, medium and heavy aircraft as the only wake-turbulence category, and once with a random subset (about 25\% of aircraft) of these artificially modeled as light aircraft so that two wake-turbulence categories are used. 
In the first hour (5:00 AM–5:59 AM), we selected 9 light aircraft, and in the second hour (6:00 AM–6:59 AM),  4 light aircraft. In addition, we consider different sizes of intervals of allowed maximum deviation from the planned time at the entry points.

\noindent{\bf Experiment Set 1.}
We divide the two hours into four half-hour intervals and  consider arrivals to the TMA between 5:00AM and 5:29AM with 16 aircraft, between 5:30AM and 5:59AM with 17 aircraft, 
between 6:00AM and 6:29AM with 9 aircraft, and between 6:30AM and 6:59AM with 7 aircraft. 
We run the experiments for the four half hours once with zero artificial light aircraft and once with six, three, two, and two simulated light aircraft in the first, second, third, and fourth half hour, respectively.
Additionally, we run these experiments with different values for the allowed maximum deviation from the planned time at the entry points.
This experiment set includes these scenarios:
\begin{itemize}
    \item[]\textbf{1a.} Run the different time periods independent of each other.
    \item[]\textbf{1b.} Include aircraft that are still in the TMA from the previous half hour, that is, aircraft that have not yet reached the runway when the next half hour begins.
    \item[]\textbf{1c.} Include tree consistency, that is, limiting the number of edges that differ between trees for consecutive time periods.
    \item[]\textbf{1d.} Include both aircraft that are still in the TMA from previous time period and tree consistency.
\end{itemize}

\noindent
{\bf Experiment Set 2.} 
We consider the two full hours in this experiment, that is, arrivals to the TMA between 5:00AM and 5:59AM with 33 aircraft and between 6:00AM and 6:59AM with 16 aircraft. 
We run the experiments for the full hours once with zero and once with  nine and four simulated light aircraft in the first and second full hour, respectively.
Additionally, we run these experiments with different values for the allowed maximum deviation from the planned time at the entry points.
Similar to Experiment 1, this experiment set consists of these scenarios:
\begin{itemize}
    \item[]\textbf{2a.} Run the different time periods independent of each other.
    \item[]\textbf{2b.} Include aircraft that are still in the TMA from the previous hour.
    \item[]\textbf{2c.} Include tree consistency. 
    \item[]\textbf{2d.} Include both aircraft that are still in the TMA from previous time period and tree consistency.
\end{itemize}

\begin{table*}[h!]
    \scriptsize
	\centering
	\caption{Results of experiment set 1a with an upper bound of 14 on the path length}\label{table:results1a-14}
	\begin{tabular}{cccccccc}\hline
                                        &                    &                    &             &                   &   Run     & Trajectory& Average \\
                                        & $|\mathcal{A}|$    &  $|C_2|$           &  $\mu$      & Tree              &  time (s) &generation & deviation \\
                                        &                    &                    &             &                   &           &time (s)   &  \\\hline
        \multirow{14}{*}{5:00--5:29}&\multirow{7}{*}{16} & \multirow{7 }{*}{6}& 0           &---                &  ---      &           &       \\[.2em]
                                        &                    &                    & 1           &$\text{T}^{a}_1$   &  137      & 18        &  0.81   \\[.2em]
                                        &                    &                    & 2           & $\text{T}^{a}_2$  &   182     & 18        & 1.19    \\[.2em]
                                        &                    &                    & 3           & $\text{T}^{a}_3$  &   164     & 18        & 1.69    \\[.2em]
                                        &                    &                    & 4           & $\text{T}^{a}_4$  &  400      &   18      &  2.12   \\[.2em]
                                        &                    &                    & 5           & $\text{T}^{a}_5$  &    107    &    18     &   2.62 \\[.2em]
                                        \cmidrule{2-8}
                                        &\multirow{7}{*}{16} & \multirow{7}{*}{0} & 0           & ---               &   ---     &           &       \\[.2em]
                                        &                    &                    & 1           & $\text{T}^{a}_6$  &   28      & 7         & 0.81    \\[.2em]
                                        &                    &                    & 2           & $\text{T}^{a}_7$  &  16       & 7         &1.37     \\[.2em]
                                        &                    &                    & 3           & $\text{T}^{a}_8$  &  13       & 7         & 2.06    \\[.2em]
                                        &                    &                    & 4           & $\text{T}^{a}_9$  &  67       & 7         &3.06     \\[.2em]
                                        &                    &                    & 5           &$\text{T}^{a}_{10}$&  93       & 7         & 3.31    \\[.2em]\hline
        \multirow{7}{*}{5:30--5:59} &\multirow{4}{*}{17} & \multirow{4}{*}{3} & $\leq 3$    & ---               &   ---     &           &       \\[.2em]
                                        &                    &                    & 4           &$\text{T}^{a}_{11}$&  435      &   22      & 2.76    \\[.2em]
                                        &                    &                    & 5           &$\text{T}^{a}_{12}$& 759       &  22       & 3.53   \\[.2em]
                                        \cmidrule{2-8}
                                        &\multirow{4}{*}{17} & \multirow{4}{*}{0} & $\leq 2$    &  ---              &  ---      &           &       \\[.2em]
                                        &                    &                    & 3           &$\text{T}^{a}_{13}$& 729       & 8         & 2.29    \\[.2em]
                                        &                    &                    & 4           &$\text{T}^{a}_{14}$& 233       &  8        & 2.88    \\[.2em]
                                        &                    &                    & 5           &$\text{T}^{a}_{15}$&  183      &  8        &  3.65  \\[.2em]\hline
        \multirow{9}{*}{6:00--6:29} &\multirow{5}{*}{9}  & \multirow{5}{*}{2}& $\leq 2$     &   ---             &  ---      &           &       \\[.2em]
                                        &                    &                    & 3           &$\text{T}^{a}_{16}$& 77        & 9         & 2.33    \\[.2em]
                                        &                    &                    & 4           &$\text{T}^{a}_{17}$&  60       & 9         & 3.22    \\[.2em]
                                        &                    &                    & 5           &$\text{T}^{a}_{18}$& 60        & 9         & 3.67   \\[.2em]
                                        \cmidrule{2-8}
                                        &\multirow{6}{*}{9}  & \multirow{6}{*}{0} &$\leq 1$     & ---               &  --       &           &       \\[.2em]
                                        &                    &                    & 2           &$\text{T}^{a}_{19}$&  18       & 4         & 1.44    \\[.2em]
                                        &                    &                    & 3           &$\text{T}^{a}_{20}$&  28       & 4         & 2.44    \\[.2em]
                                        &                    &                    & 4           &$\text{T}^{a}_{21}$&  35       &3          & 2.44    \\[.2em]
                                        &                    &                    & 5           &$\text{T}^{a}_{22}$&  27       & 3         &3.55    \\[.2em]\hline
        \multirow{14}{*}{6:30--6:59}&\multirow{7}{*}{7}  & \multirow{7}{*}{2} & 0           &   ---             &  ---      &           &       \\[.2em]
                                        &                    &                    & 1           &$\text{T}^{a}_{23}$& 5         & 7         & 0.86    \\[.2em]
                                        &                    &                    & 2           &$\text{T}^{a}_{24}$& 38        & 7         & 1.14    \\[.2em]
                                        &                    &                    & 3           &$\text{T}^{a}_{25}$& 27        & 7         & 1.71    \\[.2em]
                                        &                    &                    & 4           &$\text{T}^{a}_{26}$& 42        & 7         & 2.28    \\[.2em]
                                        &                    &                    & 5           &$\text{T}^{a}_{27}$& 47        & 7         & 4   \\[.2em]
                                        \cmidrule{2-8}
                                        &\multirow{7}{*}{7}  & \multirow{7}{*}{0} & 0           &---                &  ---      &           &       \\[.2em]
                                        &                    &                    & 1           &$\text{T}^{a}_{28}$& 3         &3          & 1    \\[.2em]
                                        &                    &                    & 2           &$\text{T}^{a}_{29}$& 7         & 3         & 1.57    \\[.2em]
                                        &                    &                    & 3           &$\text{T}^{a}_{30}$& 13        & 3         & 1.43    \\[.2em]
                                        &                    &                    & 4           &$\text{T}^{a}_{31}$&  7        & 3         & 2.57    \\[.2em]
                                        &                    &                    & 5           &$\text{T}^{a}_{32}$&  7        & 3         & 4   \\[.2em]\hline
\end{tabular}
\end{table*}

\begin{table*}[h!]
    \scriptsize
	\centering
	\caption{Results of experiment set 1b including light aircraft}\label{table:results1b-14-light}
	\begin{tabular}{ccccccccc}\hline
                                        &                    &                   &  Previous                         &          &                   &   Run     & Trajectory& Average \\
                                        & $|\mathcal{A}|$    &  $|C_2|$          & tree                              &  $\mu$   &  Tree             &  time (s) &generation & deviation \\
                                        &                    &                   &                                   &          &                   &           &time (s)   &  \\\noalign{\hrule height 1pt}
        \multirow{3}{*}{5:30--5:59} &\multirow{3}{*}{17} & \multirow{3}{*}{3}& $\text{T}^{a}_1$--$\text{T}^{a}_4$& $\leq 10$&  ---              & ---       &           &     \\[.2em]\cmidrule{4-9}
                                        &                    &                   & $\text{T}^{a}_5$                  & $\leq 8$ & ---               & ---       &           &     \\[.2em]
                                        &                    &                   & $\text{T}^{a}_5$                  & 9        &$\text{T}^{b}_1$   & 1839      &  22       & 6.06  \\[.2em]\noalign{\hrule height 1pt}
        \multirow{6}{*}{6:00--6:29} &\multirow{6}{*}{9}  & \multirow{6}{*}{2}& $\text{T}^{b}_1$                  &$\leq 2$  &---                &  ---      &           &       \\[.2em]
                                        &                    &                   &$\text{T}^{b}_1$                   & 3        &$\text{T}^{b}_2$   &   61      &    9      & 2.33   \\[.2em]
                                        &                    &                   & $\text{T}^{b}_1$                  & 4        &$\text{T}^{b}_3$   & 163       &  9        &  3.22   \\[.2em]
                                        &                    &                   & $\text{T}^{b}_1$                  & 5        &$\text{T}^{b}_4$   &  222      &    8      &  3.78  \\[.2em]
                                        &                    &                   & $\text{T}^{b}_1$                  & 6        &$\text{T}^{b}_5$   &  360      &    8      &  5.22  \\[.2em]\noalign{\hrule height 1pt}
        \multirow{17}{*}{6:30--6:59}&\multirow{17}{*}{7} &\multirow{17}{*}{2}&$\text{T}^{b}_2$                   & $\leq 4$ &---                &  ---      &           &       \\[.2em]
                                        &                    &                   &$\text{T}^{b}_2$                   & 5        &$\text{T}^{b}_6$   & 34        &7          & 3    \\[.2em]
                                        &                    &                   & $\text{T}^{b}_2$                  & 6        &$\text{T}^{b}_7$   & 64        & 7         & 4.28    \\[.2em]\cmidrule{4-9}
                                        &                    &                   & $\text{T}^{b}_3$                  & $\leq 3$ &---                & ---       &           &     \\[.2em]
                                        &                    &                   & $\text{T}^{b}_3$                  & 4        &$\text{T}^{b}_8$   & 39        & 7         & 3.43    \\[.2em]
                                        &                    &                   & $\text{T}^{b}_3$                  & 5        &$\text{T}^{b}_9$   & 91        & 7         & 4    \\[.2em]
                                        &                    &                   & $\text{T}^{b}_3$                  & 6        &$\text{T}^{b}_{10}$& 114       & 7         & 4.57    \\[.2em]\cmidrule{4-9}
                                        &                    &                   & $\text{T}^{b}_4$                  & $\leq 2$ &---                & ---       &           &     \\[.2em]
                                        &                    &                   & $\text{T}^{b}_4$                  & 3        &$\text{T}^{b}_{11}$& 34        & 7         & 2.14    \\[.2em]
                                        &                    &                   & $\text{T}^{b}_4$                  & 4        &$\text{T}^{b}_{12}$& 43        & 7         & 3.14    \\[.2em]
                                        &                    &                   & $\text{T}^{b}_4$                  & 5        &$\text{T}^{b}_{13}$& 62        & 7         & 4.14   \\[.2em]
                                        &                    &                   & $\text{T}^{b}_4$                  & 6        &$\text{T}^{b}_{14}$& 69        & 7         & 5.28   \\[.2em]\cmidrule{4-9}
                                        &                    &                   & $\text{T}^{b}_5$                  & $\leq 4$ &---                & ---       &           &  \\[.2em]
                                        &                    &                   & $\text{T}^{b}_5$                  & 5        &$\text{T}^{b}_{15}$& 61        & 7         & 3.43   \\[.2em]
                                        &                    &                   & $\text{T}^{b}_5$                  & 6        &$\text{T}^{b}_{16}$& 94        & 7         & 4   \\ [.2em]\hline
\end{tabular}
\end{table*}

The computational results of all experiments with one category $C_1$ and two categories $C_1$ and $C_2$ as well as with different values for the entry point time-window offset $\mu$ are summarized in Tables~\ref{table:results1a-14}--\ref{table:results2d-14-light} and Tables~\ref{table:results1b-14}--\ref{table:results2d-14} in Appendix~\ref{appA}.  
Results for Experiments 1a and 2a appear in Tables~\ref{table:results1a-14} and~\ref{table:results2a-14}, respectively; both tables include the single-category and two-category cases. 
As the single-category experiments are computationally easier than the two-category cases, 
we report only the harder scenarios here and include the complete tables for the single-category cases in Appendix~\ref{appA}.
Experiments 1b--1d and 2b--2d with two categories are reported in Tables~\ref{table:results1b-14-light}--\ref{table:results1d-14-light} and~\ref{table:results2b-14-light}--\ref{table:results2d-14-light}, 
while the corresponding single-category results are given in Tables~\ref{table:results1b-14}--\ref{table:results1d-14} and~\ref{table:results2b-14}--\ref{table:results2d-14} in Appendix~\ref{appA}.
In these tables, the total number of aircraft and the number of light aircraft are denoted by  $|\mathcal{A}|$ and $|C_2|$, respectively. 
(Where $|C_2|=0$ indicates that we have only one aircraft category.)
Moreover, the time required to solve the model is denoted as ``run time''  while time needed to generate the sets $\Xi$ and $\Upsilon$ is given as ``trajectory generation time''. 
Both times are given in seconds.
The column labeled ``tree'' shows the name of the optimal tree for each row, where ``---'' indicates that the problem for that row is infeasible.
We label the trees for scenarios a, b, c, and d as $T^a$, $T^b$, $T^c$, and $T^d$, respectively.
Furthermore, when additional constraints (such as tree consistency and the presence of aircraft carried over from the previous time period in the TMA) are included, 
the table contains at most two extra columns: ``previous tree'', which gives the name of the previous tree, and $U$, that is, the maximum number of edges that may differ between the two trees.
The final column in each table shows the average entry-time deviation, defined as the mean 
absolute value of the difference between the original time of aircraft arrival to the TMA and the scheduled entry time according to our optimized arrival schedule.

\begin{table*}[h!]
    \scriptsize
	\centering
		\setlength{\tabcolsep}{5pt}
	\caption{Results of experiment set 1c including light aircraft}\label{table:results1c-14-light}
	\begin{tabular}{cccccccccc}\hline
                                        &                    &                      &   Previous        &          &            &                   &       Run    &  Trajectory & Average \\
                                        & $|\mathcal{A}|$    &  $|C_2|$             & tree              & $\mu$    &  $U$       & Tree              &   time (s)   &generation   & deviation \\
                                        &                    &                      &                   &          &            &                   &              &time (s)     &  \\\noalign{\hrule height 1pt}
        \multirow{26}{*}{5:30--5:59}&\multirow{26}{*}{17}& \multirow{26}{*}{3}  & $\text{T}^{a}_1$  & $\leq 5$ & $\leq 10$  & ---               &   ---        &             &     \\[.2em]
                                        &                    &                      & $\text{T}^{a}_1$  & 6        &$\leq 9$    &---                &  ---         &             &     \\[.2em]
                                        &                    &                      & $\text{T}^{a}_1$  & 6        &10          &$\text{T}^{c}_1$   & 384          & 22          &  4.59 \\[.2em]
                                        &                    &                      & $\text{T}^{a}_1$  & 7        & $\leq 4$   &---                &  ---         &             &     \\[.2em]
                                        &                    &                      & $\text{T}^{a}_1$  & 7        &   5        &$\text{T}^{c}_2$   & 318          & 27          &  5 \\[.2em]\cmidrule{4-10}
                                        &                    &                      & $\text{T}^{a}_2$  & $\leq 3$ & $\leq 10$  &---                &  ---         &             &     \\[.2em]
                                        &                    &                      & $\text{T}^{a}_2$  & 4        & $\leq 7$   &---                &  ---         &             &     \\[.2em]
                                        &                    &                      & $\text{T}^{a}_2$  & 4        & 8          &$\text{T}^{c}_3$   &  513         & 22          &  2.41   \\[.2em]
                                        &                    &                      & $\text{T}^{a}_2$  & 5        & $\leq 3$   &---                &  ---         &             &     \\[.2em]
                                        &                    &                      & $\text{T}^{a}_2$  & 5        & 4          &$\text{T}^{c}_4$   & 478          & 23          &  3.82   \\[.2em]
                                        &                    &                      & $\text{T}^{a}_2$  & 6        & 1          &$\text{T}^{c}_5$   & 132          &  22         &  4.65 \\[.2em]\cmidrule{4-10}
                                        &                    &                      & $\text{T}^{a}_3$  & $\leq 3$ & $\leq 10$  &---                &  ---         &             &     \\[.2em]
                                        &                    &                      & $\text{T}^{a}_3$  & 4        & $\leq 7$   &---                &  ---         &             &     \\[.2em]
                                        &                    &                      & $\text{T}^{a}_3$  & 4        & 8          &$\text{T}^{c}_6$   &  387         &  22         & 2.82    \\[.2em]
                                        &                    &                      & $\text{T}^{a}_3$  & 5        & $\leq 3$   &---                &  ---         &             &     \\[.2em]
                                        &                    &                      & $\text{T}^{a}_3$  & 5        & 4          &$\text{T}^{c}_7$   &  468         &   22        &  4.18   \\[.2em]\cmidrule{4-10}
                                        &                    &                      & $\text{T}^{a}_4$  &$\leq 3$  & $\leq 10$  &---                &  ---         &             &     \\[.2em]
                                        &                    &                      & $\text{T}^{a}_4$  & 4        & $\leq 9$   &---                &  ---         &             &     \\[.2em]
                                        &                    &                      & $\text{T}^{a}_4$  & 4        & 10         &$\text{T}^{c}_8$   &  412         & 22          & 2.94    \\[.2em]
                                        &                    &                      & $\text{T}^{a}_4$  & 5        & $\leq 7$   &---                &  ---         &             &     \\[.2em]
                                        &                    &                      & $\text{T}^{a}_4$  & 5        &8           &$\text{T}^{c}_9$   &  308         & 22          & 3.65\\[.2em]\cmidrule{4-10}
                                        &                    &                      & $\text{T}^{a}_5$  & $\leq 4$ & $\leq 10$  & ---               &  ---         &             &     \\[.2em]
                                        &                    &                      & $\text{T}^{a}_5$  & 5        & 1          &$\text{T}^{c}_{10}$&  107         & 22          &  3.59   \\[.2em]
                                                                                    \noalign{\hrule height 1pt}
        \multirow{8}{*}{6:00--6:29} &\multirow{8}{*}{9}  & \multirow{8}{*}{2}   &$\text{T}^{c}_{10}$& 3        &  $\leq 15$ & ---               &  ---         &             &       \\[.2em]
                                        &                    &                      &$\text{T}^{c}_{10}$& 4        &  $\leq 10$ &  ---              &   ---        &             &     \\[.2em]
                                        &                    &                      &$\text{T}^{c}_{10}$& 4        &   11       &$\text{T}^{c}_{11}$& 137          &   9         &   3.44  \\[.2em]
                                        &                    &                      &$\text{T}^{c}_{10}$& 5        &  $\leq 10$ &  ---              &   ---        &             &     \\[.2em]
                                        &                    &                      &$\text{T}^{c}_{10}$& 5        &   11       &$\text{T}^{c}_{12}$& 132          &   9         &   3.77  \\[.2em]
                                        &                    &                      &$\text{T}^{c}_{10}$& 6        &  $\leq 3$  &  ---              &   ---        &             &     \\[.2em]
                                        &                    &                      &$\text{T}^{c}_{10}$& 6        &  4         &$\text{T}^{c}_{13}$& 152          &   9         &   5.33  \\[.2em]\noalign{\hrule height 1pt}
        \multirow{6}{*}{6:30--6:59} &\multirow{6}{*}{7}  & \multirow{6}{*}{2}   &$\text{T}^{c}_{11}$& 1        & 1          &$\text{T}^{c}_{14}$&  4           &  6          &   1    \\[.2em]
                                        &                    &                      &$\text{T}^{c}_{12}$& 1        & 1          &$\text{T}^{c}_{15}$&  4           & 6           & 0.86    \\[.2em]
                                        &                    &                      &$\text{T}^{c}_{13}$& 1        & $\leq 3$   & ---               & ---          &             &        \\[.2em]
                                        &                    &                      &$\text{T}^{c}_{13}$& 1        & 4          &$\text{T}^{c}_{16}$& 5            & 6           & 1    \\[.2em]
                                        &                    &                      &$\text{T}^{c}_{13}$& 2        & 1          &$\text{T}^{c}_{17}$& 7            & 6           & 1.57    \\[.2em]
                                       \hline
\end{tabular}
\end{table*}

We first ran Experiment~1a with two sets of paths: one using an upper bound of 14 edges and one using an upper bound of 15 edges. 
Table~\ref{table:results1a-14} and~\ref{table:results1a-15} in Appendix~\ref{appA} show the results of Experminet~1a based on  paths with a maximum of 14 and 15 edges, respectively.
The optimal values of the first case are obtained within 5 to 759 seconds while for the second case they are obtained within 26 to 3129 seconds.
Although both variants produce identical objective values, the 14-edge variant is consistently much faster, 
which shows that computational time  highly depends on the path-length bound. 
Therefore, we used a maximum of 14 edges for all subsequent experiments to reduce runtime without affecting solution quality.

\begin{table*}[h!]
    \scriptsize
	\centering
		\setlength{\tabcolsep}{5pt}
	\caption{Results of experiment set 1d including light aircraft}\label{table:results1d-14-light}
	\begin{tabular}{cccccccccc}\hline
                                        &                     &                         &   Previous       &          &            &                 &       Run    &  Trajectory & Average \\
                                        & $|\mathcal{A}|$     &  $|C_2|$                & tree             & $\mu$    &  $U$       & Tree            &   time (s)   &generation   & deviation \\
                                        &                     &                         &                  &          &            &                 &              &time (s)     &  \\[.2em]\noalign{\hrule height 1pt}
        \multirow{4}{*}{5:30--5:59} &\multirow{4}{*}{17}  & \multirow{4}{*}{3}      & $\text{T}^{a}_5$ & $\leq 8$ & $\leq 15$  & ---             &   ---        &             &     \\[.2em]
                                        &                     &                         & $\text{T}^{a}_5$ & 9        & $\leq 8$   & ---             &  ---         &             &     \\[.2em]
                                        &                     &                         & $\text{T}^{a}_5$ & 9        & 9          &$\text{T}^{d}_1$ &  1596        & 22          &  6.53   \\[.2em]
                                \noalign{\hrule height 1pt}
        \multirow{4}{*}{6:00--6:29} &\multirow{4}{*}{9}   & \multirow{4}{*}{2}      &$\text{T}^{d}_1$  &$\leq 6$  & $\leq 15$  & ---             &   ---        &             &     \\[.2em]
                                        &                     &                         &$\text{T}^{d}_1$  & 7        &  $\leq 2$  &  ---            & ---          &             &     \\[.2em]
                                        &                     &                         &$\text{T}^{d}_1$  & 7        &    3       &$\text{T}^{d}_2$ & 90           &   9         &  5.11  \\[.2em]
                            \noalign{\hrule height 1pt}
        \multirow{4}{*}{6:30--6:59} &\multirow{4}{*}{7}   & \multirow{4}{*}{2}      &$\text{T}^{d}_2$  &$\leq 5$  & $\leq 15$  & ---             &   ---        &             &     \\[.2em]
                                        &                     &                         &$\text{T}^{d}_2$  & 6        &  $\leq 8$  &  ---            & ---          &             &     \\[.2em]
                                        &                     &                         &$\text{T}^{d}_2$  & 6        &    9       &$\text{T}^{d}_3$ & 121          &   7         &  3.71  \\[.2em]
                                        \hline
\end{tabular}
\end{table*}

\begin{table*}[h!]
    \scriptsize
	\centering
	\caption{Results of experiment set 2a}\label{table:results2a-14}
	\begin{tabular}{cccccccc}\hline
                                        &                    &                      &           &                    &    Run     &  Trajectory  & Average \\
                                        & $|\mathcal{A}|$    &  $|C_2|$             &  $\mu$    & Tree               &   time (s) &generation    & deviation \\
                                        &                    &                      &           &                    &            &time (s)      &  \\[.2em]\hline
        \multirow{5}{*}{5:00--5:59} &\multirow{2}{*}{33} & \multirow{2}{*}{9}   & $\leq 7$  &  ---               &  ---       &              &       \\[.2em]
                                        &                    &                      & 8         &$\text{T}^{a}_{66}$ &     3514   &   75         &   5.51     \\[.2em]\cmidrule{2-8}
                                        &\multirow{3}{*}{33} & \multirow{3}{*}{0}   & $\leq 5$  &  ---               &   ---      &              &       \\[.2em]
                                        &                    &                      & 6         &$\text{T}^{a}_{67}$ &  2875      &   25         &   4.33     \\[.2em]
                                        &                    &                      & 7         &$\text{T}^{a}_{68}$ &  1333      &   25         &   4.15     \\[.2em]\hline
        \multirow{7}{*}{6:00--6:59} &\multirow{4}{*}{16} & \multirow{4}{*}{4}   & $\leq 3$  &  ---               &  ---       &              &       \\[.2em]
                                        &                    &                      & 4         &$\text{T}^{a}_{69}$ &  255       & 30           & 3.12    \\[.2em]
                                        &                    &                      & 5         &$\text{T}^{a}_{70}$ & 357        & 30           & 3.25   \\[.2em]\cmidrule{2-8}
                                        &\multirow{6}{*}{16} & \multirow{6}{*}{0}   & $\leq 2$  &  ---               &  ---       &              &       \\[.2em]
                                        &                    &                      & 3         &$\text{T}^{a}_{71}$ &  52        &   12         &   2.12   \\[.2em]
                                        &                    &                      & 4         &$\text{T}^{a}_{72}$ &  27        & 12           & 3.5    \\[.2em]
                                        &                    &                      & 5         &$\text{T}^{a}_{73}$ &  36        &  12          &3.25   \\[.2em]
                                        &                    &                      & 6         &$\text{T}^{a}_{74}$ &  37        &  12          &4.67   \\[.2em]\hline
\end{tabular}
\end{table*}

\begin{table*}[h!]
    \scriptsize
	\centering
	\caption{Results of experiment set 2b including light aircraft}\label{table:results2b-14-light}
	\begin{tabular}{ccccccccc}\hline
                                        &                    &                      &  Previous         &            &                    &     Run     &  Trajectory   &  Average\\
                                        & $|\mathcal{A}|$    &  $|C_2|$             & tree              &  $\mu$     & Tree               &   time (s)  &generation     & deviation \\
                                        &                    &                      &                   &            &                    &             &time (s)       &  \\\hline
        \multirow{4}{*}{6:00--6:59}&\multirow{4}{*}{16} & \multirow{4}{*}{4}  &$\text{T}^{a}_{66}$& $\leq 3$   & ---                & ---         &               &       \\
                                        &                    &                      &$\text{T}^{a}_{66}$& 4          &$\text{T}^{b}_{26}$ &  533        &  26           &  3.06      \\
                                        &                    &                      &$\text{T}^{a}_{66}$& 5          &$\text{T}^{b}_{27}$ &  909        &  26           &  3.56       \\
                                        &                    &                      &$\text{T}^{a}_{66}$& 6          &$\text{T}^{b}_{28}$ &  722        &  26           &  4.25       \\\hline
    \end{tabular}
\end{table*}

\begin{table*}[h!]
    \scriptsize
	\centering
	\setlength{\tabcolsep}{5pt}
	\caption{Results of experiment set 2c including light aircraft}\label{table:results2c-14-light}
	\begin{tabular}{cccccccccc}\hline
                                        &                    &                      & Previous          &           &           &                    &     Run   & Trajectory&  Average\\
                                        & $|\mathcal{A}|$    &  $|C_2|$             & tree              &  $\mu$    &$U$        & Tree               &   time (s)&generation & deviation \\
                                        &                    &                      &                   &           &           &                    &           &time (s)   &  \\\noalign{\hrule height 1pt}
        \multirow{5}{*}{6:00--6:59}&\multirow{5}{*}{16} & \multirow{5}{*}{4}   &$\text{T}^{a}_{66}$& $\leq 4$  &$\leq 12$  & ---                &  ---      &           &       \\
                                        &                    &                      &$\text{T}^{a}_{66}$& 5         &$\leq 10$  & ---                &  ---      &           &        \\
                                        &                    &                      &$\text{T}^{a}_{66}$& 5         &11         &$\text{T}^{c}_{40}$ & 240       &  27       &  3.5      \\
                                        &                    &                      &$\text{T}^{a}_{66}$& 6         &$\leq 3$   & ---                &  ---      &           &        \\
                                        &                    &                      &$\text{T}^{a}_{66}$& 6         &4          &$\text{T}^{c}_{41}$ &  288      &   27      &    4.62    \\
                                        \hline
\end{tabular}
\end{table*}

\begin{table*}[h!]
    \scriptsize
	\centering
		\setlength{\tabcolsep}{5pt}
	\caption{Results of experiment set 2d including light aircraft}\label{table:results2d-14-light}
	\begin{tabular}{cccccccccc}\hline                                        
                                        &                    &                      &  Previous         &           &           &                    &       Run & Trajectory&  Average\\
                                        & $|\mathcal{A}|$    &  $|C_2|$             & tree              &  $\mu$    &$U$        & Tree               &   time (s)&generation & deviation \\
                                        &                    &                      &                   &           &           &                    &           &time (s)   &  \\\hline
        \multirow{5}{*}{6:00--6:59}&\multirow{5}{*}{16} & \multirow{5 }{*}{4}  &$\text{T}^{a}_{66}$& $\leq 5$  &$\leq 15$  & ---                &  ---      &           &       \\
                                        &                    &                      &$\text{T}^{a}_{66}$& 6         &$\leq 14$  & ---                &  ---      &           &        \\
                                        &                    &                      &$\text{T}^{a}_{66}$& 6         &15         &$\text{T}^{d}_7$    &  678      &  27       &  4.25      \\
                                        &                    &                      &$\text{T}^{a}_{66}$& 7         &$\leq 4$   & ---                &  ---      &           &        \\
                                        &                    &                      &$\text{T}^{a}_{66}$& 7         &5          &$\text{T}^{d}_8$    &  469      &   27      &    6.06    \\
                                        \hline
\end{tabular}
\end{table*}

In the consecutive high-traffic scenarios, feasibility becomes challenging when aircraft remaining in the TMA from the previous period must be included.
For instance, the results of Experiment~1b (in the two-category setting)  presented in Table~\ref{table:results1b-14-light} show that 
a feasible routing tree for the time period 5:30--5:59 was obtained only with an entry point offset of 9 minutes,
whereas the preceding half-hour used  an entry point offset of 5 minutes.
Overall, the results  indicate that each additional operational feature increases the minimum value of $\mu$ required for feasibility.

\begin{figure*}[t] \hspace*{-1cm}
	\begin{subfigure}[t]{0.35\textwidth}
	\centering
	\includegraphics[width=0.5\textwidth]{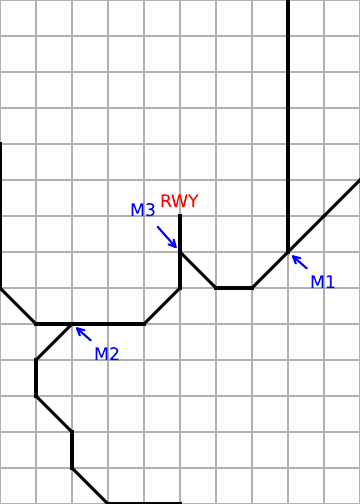} \caption{$\mu=1$}\label{fig:1a-1}
\end{subfigure}\hfill \hspace*{-2.6cm}
\begin{subfigure}[t]{0.35\textwidth}
	\centering
	\includegraphics[width=.5\textwidth]{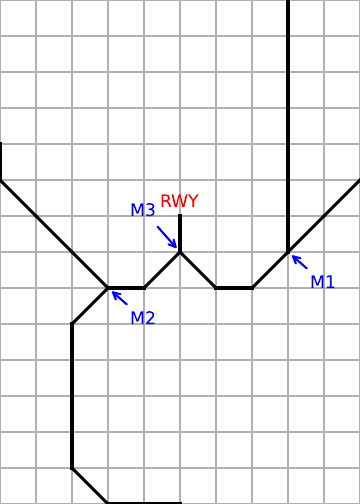} \caption{$\mu=2$}\label{fig:1a-2}
\end{subfigure}\hfill \hspace*{-2.6cm}
\begin{subfigure}[t]{0.35\textwidth}
	\centering
	\includegraphics[width=.5\textwidth]{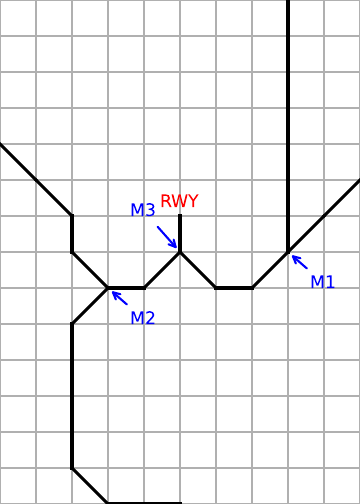} \caption{$\mu=3$}\label{fig:1a-3}
\end{subfigure}\hfill \hspace*{-2.6cm}
\begin{subfigure}[t]{0.35\textwidth}
	\centering
	\includegraphics[width=.5\textwidth]{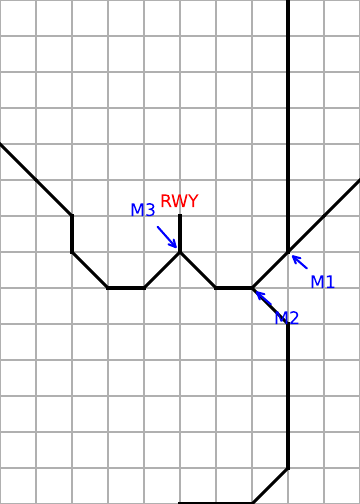} \caption{$\mu=4$}\label{fig:1a-4}
\end{subfigure}\hfill\hspace*{-2.6cm}
\begin{subfigure}[t]{0.35\textwidth}
	\centering
	\includegraphics[width=.5\textwidth]{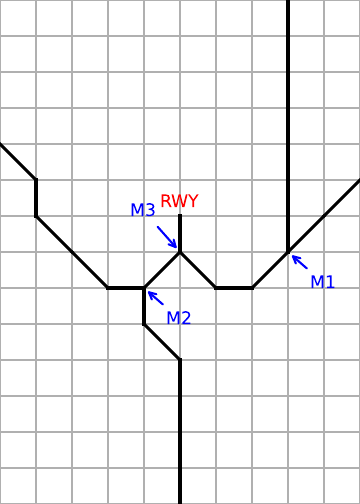} \caption{$\mu=5$}\label{fig:1a-5}
\end{subfigure}\hspace*{-1.25cm}
\caption{Trees of Experminet~1a for period 5:00--5:29 including light aircraft}\label{fig:1a-500-529}
\end{figure*}

\begin{figure*}[t] \hspace*{-.5cm}
	\begin{subfigure}[t]{0.35\textwidth}
	\centering
	\includegraphics[width=0.5\textwidth]{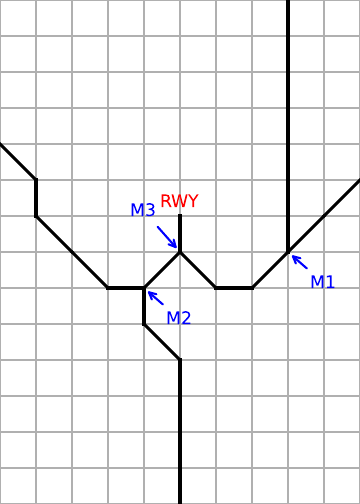} \caption{5:00--5:29}\label{fig:1a1-5}
\end{subfigure}\hfill \hspace*{-2.6cm}
\begin{subfigure}[t]{0.35\textwidth}
	\centering
	\includegraphics[width=.5\textwidth]{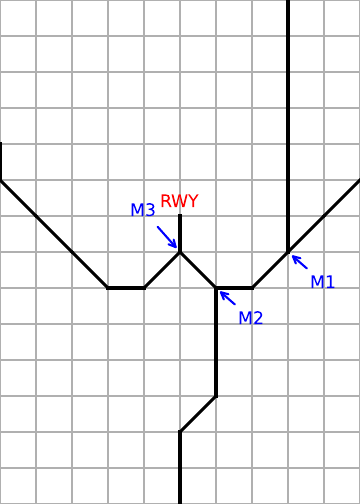} \caption{5:30--5:59}\label{fig:1a2-5}
\end{subfigure}\hfill \hspace*{-2.6cm}
\begin{subfigure}[t]{0.35\textwidth}
	\centering
	\includegraphics[width=.5\textwidth]{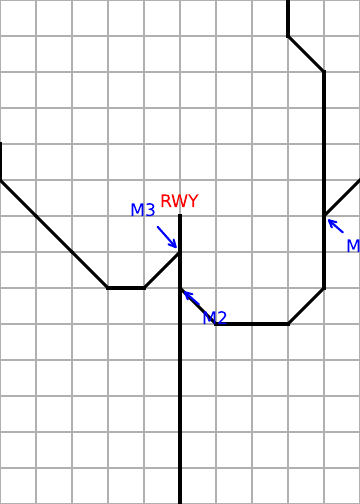} \caption{6:00--6:29}\label{fig:1a3-5}
\end{subfigure}\hfill \hspace*{-2.6cm}
\begin{subfigure}[t]{0.35\textwidth}
	\centering
	\includegraphics[width=.5\textwidth]{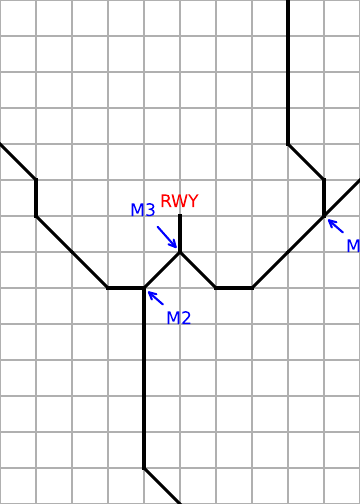} \caption{6:30--6:59}\label{fig:1a4-5}
\end{subfigure}\hspace*{-.5cm}
\caption{Trees of Experiment~1a including light aircraft with $\mu=5$}\label{fig:1a-500-600}
\end{figure*}

\begin{figure*}[t] \hspace*{-.5cm}
	\begin{subfigure}[t]{0.35\textwidth}
	\centering
	\includegraphics[width=0.5\textwidth]{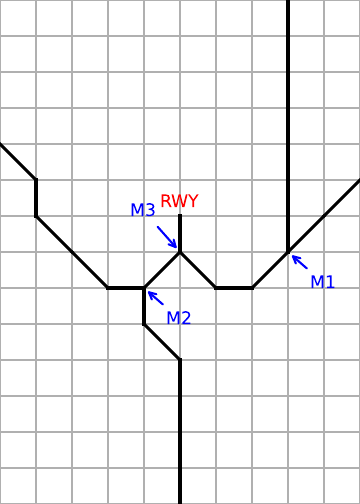} \caption{Tree $\text{T}_5^a$}\label{fig:1b-5}
\end{subfigure}\hfill \hspace*{-2.6cm}
\begin{subfigure}[t]{0.35\textwidth}
	\centering
	\includegraphics[width=.5\textwidth]{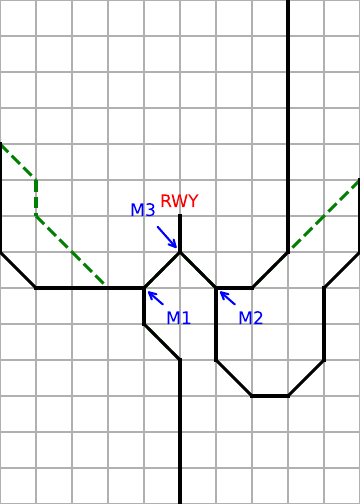} \caption{Tree $\text{T}_1^b$}\label{fig:1b-59}
\end{subfigure}\hfill \hspace*{-2.6cm}
\begin{subfigure}[t]{0.35\textwidth}
	\centering
	\includegraphics[width=.5\textwidth]{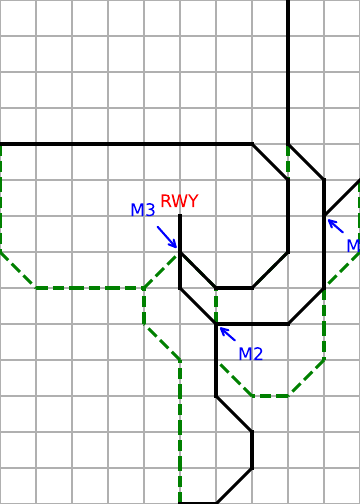} \caption{Tree $\text{T}_5^b$}\label{fig:1b-595}
\end{subfigure}\hfill \hspace*{-2.6cm}
\begin{subfigure}[t]{0.35\textwidth}
	\centering
	\includegraphics[width=.5\textwidth]{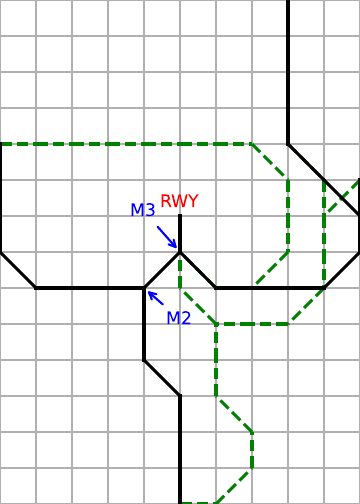} \caption{Tree $\text{T}_{15}^a$}\label{fig:1b-5955}
\end{subfigure}\hspace*{-.5cm}
\caption{Some trees of Experiment~1b including light aircraft (the tree from the previous time period is indicated in green dashed lines)}\label{fig:1b}
\end{figure*}

\begin{figure*}[t] \hspace*{-.5cm}
	\begin{subfigure}[t]{0.35\textwidth}
	\centering
	\includegraphics[width=0.5\textwidth]{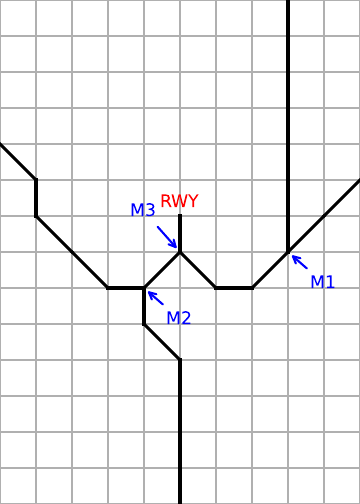}\caption{Tree $\text{T}_5^a$}\label{fig:1c-5}
\end{subfigure}\hfill \hspace*{-2.6cm}
\begin{subfigure}[t]{0.35\textwidth}
	\centering
	\includegraphics[width=.5\textwidth]{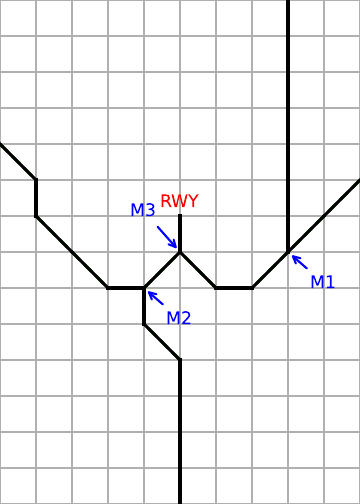}\caption{Tree $\text{T}_{10}^c$}\label{fig:1c-55}
\end{subfigure}\hfill \hspace*{-2.6cm}
\begin{subfigure}[t]{0.35\textwidth}
	\centering
	\includegraphics[width=.5\textwidth]{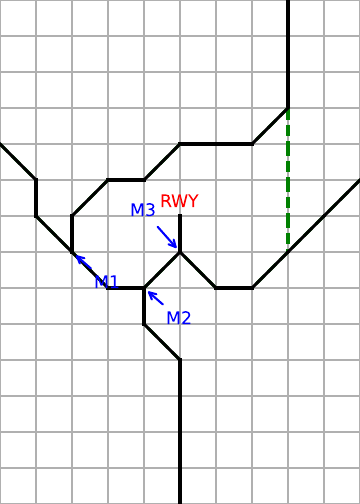} \caption{Tree $\text{T}_{12}^c$}\label{fig:1c-555}
\end{subfigure}\hfill \hspace*{-2.6cm}
\begin{subfigure}[t]{0.35\textwidth}
	\centering
	\includegraphics[width=.5\textwidth]{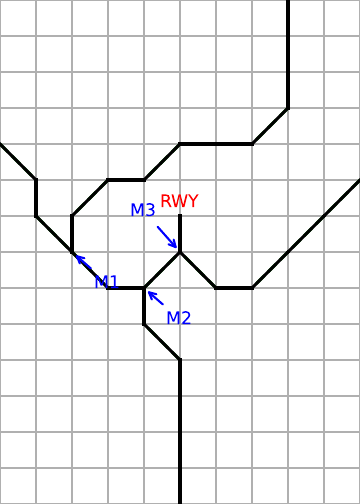} \caption{Tree $\text{T}_{15}^c$}\label{fig:1c-5551}
\end{subfigure}\hspace*{-.5cm}
\caption{Some trees of Experiment~1c including light aircraft}\label{fig:1c}
\end{figure*}

\begin{figure*}[t] \hspace*{-.5cm}
	\begin{subfigure}[t]{0.35\textwidth}
	\centering
	\includegraphics[width=0.5\textwidth]{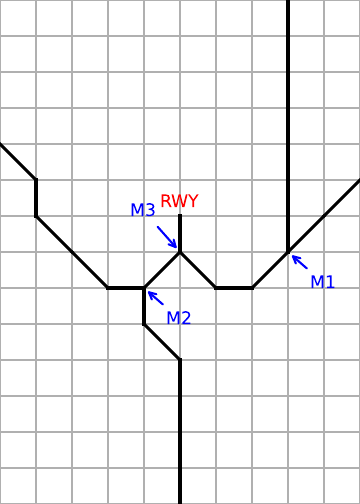} \caption{Tree $\text{T}_5^a$}\label{fig:1d-5}
\end{subfigure}\hfill \hspace*{-2.6cm}
\begin{subfigure}[t]{0.35\textwidth}
	\centering
	\includegraphics[width=.5\textwidth]{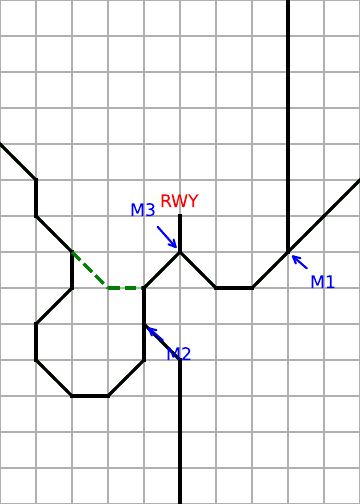} \caption{Tree $\text{T}_1^d$}\label{fig:1d-59}
\end{subfigure}\hfill \hspace*{-2.6cm}
\begin{subfigure}[t]{0.35\textwidth}
	\centering
	\includegraphics[width=.5\textwidth]{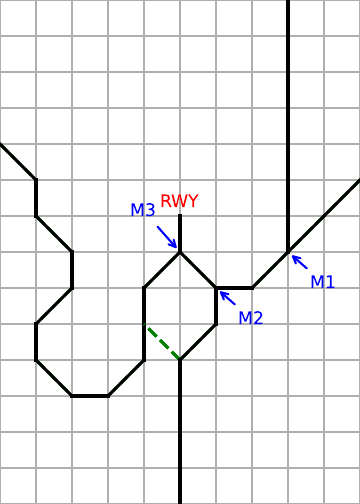} \caption{Tree $\text{T}_2^d$}\label{fig:1d-597}
\end{subfigure}\hfill \hspace*{-2.6cm}
\begin{subfigure}[t]{0.35\textwidth}
	\centering
	\includegraphics[width=.5\textwidth]{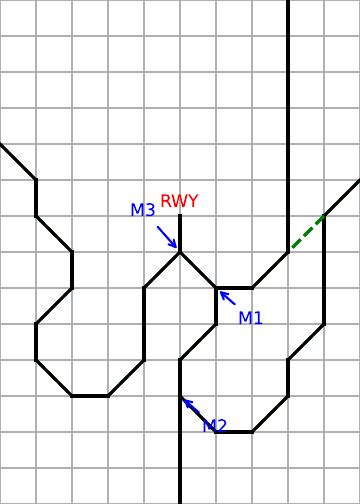} \caption{Tree $\text{T}_3^d$}\label{fig:1d-5976}
\end{subfigure}\hspace*{-.5cm}
\caption{Some trees of Experiment~1d including light aircraft}\label{fig:1d}
\end{figure*}

Figure~\ref{fig:1a-500-529} shows the arrival trees from Experiment~1a for time period 5:00--5:29 including light aircraft with different maximum deviation from the planned time at the entry points. 
The figure illustrates  how allowing larger deviations from the planned entry time yields progressively more compact and regular trees. 
In Figure~\ref{fig:1a-500-600}, the trees from Experiment~1a for each half-hour period  are shown for $\mu=5$, 
highlighting the stability of the tree structure when they are independently solved  with the same entry time window. 
Figure~\ref{fig:1b} presents the sequence of trees obtained in Experiment~1b including light aircraft. 
The tree in Figure~\ref{fig:1b-5} is obtained for time period 5:00--5:29 and forms  the basis for Figure~\ref{fig:1b-59}.
Then, the tree in Figure~\ref{fig:1b-59} which is for time period 5:30--5:59 took into account aircraft presence in the TMA from the previous time period, that is, from the tree in Figure~\ref{fig:1b-5}. 
The previous tree (i.e., Figure~\ref{fig:1b-5}) is overlaid on Figure~\ref{fig:1b-59} as a green dashed line.
Similarly, the tree in Figure~\ref{fig:1b-595}, which is for time period 6:00--6:29, considered the  presence of aircraft from the preceding period Figure~\ref{fig:1b-59}.
Then, the tree in Figure~\ref{fig:1b-5955} is for the time period 6:30--6:59 and is influenced by the result shown in Figure~\ref{fig:1b-595}.
In each figure, the previous tree is overlaid as a green dashed line.
A comparison between Figure~\ref{fig:1b} and Figure~\ref{fig:1a-500-600} highlights the significant influence of aircraft remaining from the previous period on route selection. 
Furthermore, Figures~\ref{fig:1c} and~\ref{fig:1d} illustrate the trees for Experiments~1c and~1d, respectively. 
These figures show that a larger $\mu$ is required in order to get consistency in trees across successive periods.

\vspace{1cm}

Finally, Table~\ref{table:results-T5a} presents the optimized time schedule for the DARs 
obtained in Experiment~1a during the 05:00--05:29 period, including light aircraft with $\mu=5$.
Light aircraft are indicated by the color red in the table.
In addition, entry points HMR (north), ELTOK (west), XILAN (east), and NILUG (south) are labeled by Ent1, Ent2, Ent3, and Ent4,  respectively. 
For each aircraft, the table lists the time of entry at the designated entry point, the time at which it passes each merge point, and its final arrival time at runway.
The table confirms that  all temporal separation  constraints are satisfied along every path.
Similary, in Table~\ref{table:results-T1b}, we show the optimized time schedule for the  tree obtained in Experiment~1b during the 5:30--5:59 period,
which considers aircraft remaining in the TMA from the previous period. 
A comparison between  Table~\ref{table:results-T1b} and Table~\ref{table:results-T5a} shows that the temporal separation is fulfilled even when transitioning from one  tree to the next.

\begin{table}[h!]
     \small
    	\centering
    	    \caption{Optimized time schedule for the tree $\text{T}_5^a$}\label{table:results-T5a}
    \begin{tabular}{lccccccc}\hline
        Entry& Aircraft & $\bar{t}_a$ & Scheduled & M1 & M2 & M3& RWY\\
        point & &  & time &  &  & & \\\hline
        Ent1 & 6 & 5:03 & 4:59 & 5:06 & -- & 5:09 & 5:10\\
        Ent1 & 8 & 5:12 & 5:17 & 5:24 & -- & 5:27 & 5:29\\
        Ent1 & 10 & 5:23 & 5:24 & 5:31 & -- & 5:34 & 5:35\\
        \rowcolor{red!15} Ent1 & 32 & 5:06 & 5:06 & 5:13 & -- & 5:16 & 5:18\\
        \rowcolor{red!15} Ent1 & 34 & 5:09 & 5:08 & 5:15 & -- & 5:18 & 5:20\\
        Ent2 & 14 & 5:02 & 4:57 & -- & 5:02 & 5:03 & 5:04\\
        Ent2 & 16 & 5:07 & 5:07 & -- & 5:12 & 5:13 & 5:14\\
        Ent2 & 18 & 5:22 & 5:25 & -- & 5:30 & 5:31 & 5:33\\
        \rowcolor{red!15} Ent2 & 35 & 5:26 & 5:31 & -- & 5:36 & 5:37 & 5:39\\
        \rowcolor{red!15} Ent2 & 37 & 5:15 & 5:14 & -- & 5:19 & 5:20 & 5:22\\
        Ent3 & 3 & 5:21 & 5:24 & 5:26 & -- & 5:29 & 5:31\\
        \rowcolor{red!15} Ent3 & 27 & 5:29 & 5:34 & 5:36 & -- & 5:39 & 5:41\\
        \rowcolor{red!15} Ent3 & 28 & 5:05 & 5:01 & 5:03 & --& 5:06 & 5:08\\
        Ent4 & 23 & 5:17 & 5:18 & -- & 5:24 & 5:25 & 5:26\\
        Ent4 & 24 & 5:06 & 5:04 & -- & 5:10 & 5:11 & 5:12\\
        Ent4 & 25 & 5:14 & 5:16 & -- & 5:22 & 5:23 & 5:24\\\hline
    \end{tabular}
\end{table}

\begin{table}[h!]
     \small
    \centering
      \caption{Optimized time schedule for the tree $\text{T}_1^b$}\label{table:results-T1b}
    \begin{tabular}{lccccccc}\hline
        Entry& Aircraft & $\bar{t}_a$ & Scheduled & M1 & M2 & M3& RWY\\
        point & &  & time &  &  & & \\\hline
        Ent1 & 12 & 5:41 & 5:38 & -- & 5:47 & 5:49 & 5:51\\
        Ent1 & 13 & 5:40 & 5:48 & -- & 5:57 & 5:58 & 5:59\\
        \rowcolor{red!15} Ent1 & 30 & 5:37 & 5:31 & -- & 5:40 & 5:41 & 5:43\\
        Ent1 & 31 & 5:34 & 5:35 & -- & 5:44 & 5:45 & 5:47\\
        Ent1 & 61 & 5:47 & 5:56 & -- & 6:05 & 6:06 & 6:07\\
        Ent1 & 63 & 5:52 & 5:59 & -- & 6:08 & 6:10 & 6:12\\
        Ent1 & 65 & 5:44 & 5:52 & -- & 6:01 & 6:02 & 6:03\\
        Ent2 & 19 & 5:36 & 5:43 & 5:51 & -- & 5:52 & 5:53\\
        Ent2 & 20 & 5:31 & 5:39 & 5:46 & -- & 5:47 & 5:49\\
        Ent2 & 50 & 5:45 & 5:46 & 5:53 & -- & 5:54 & 5:55\\
        Ent2 & 52 & 5:47 & 5:56 & 6:03 & -- & 6:04 & 6:05\\
        Ent2 & 55 & 5:44 & 5:52 & 5:59 & -- & 6:00 & 6:01\\
     \rowcolor{red!15} Ent2 & 69 & 5:57 & 6:05 & 6:12 & -- & 6:13 & 6:15\\
        Ent3 & 5 & 5:33 & 5:42 & -- & 6:12 & 6:15 & 6:18\\
        \rowcolor{red!15} Ent4 & 39 & 5:39 & 5:36 & 5:42 & -- & 5:43 & 5:45\\
        Ent4 & 46 & 5:45 & 5:49 & 5:55 & -- & 5:56 & 5:57\\
        Ent4 & 47 & 5:57 & 6:01 & 6:07 & -- & 6:08 & 6:09\\\hline
  \end{tabular}
\end{table}

In comparison with the results presented by~\cite{spshspp-asmda-21}, our model significantly outperforms the previous model in terms of computational efficiency. 
The largest instances that the model by~\cite{spshspp-asmda-21} could solve were each half-hour case in Experiment~1a with a maximum of 14 edges. 
We managed to solve those instances within 5 seconds to 12.65 minutes, while the previous model required 40.9 hours for each. 
We should also mention that only 30 out of the 33 aircraft have been considered by~\cite{spshspp-asmda-21}. 
In contrast, our model can handle the complete traffic for the full hour, including all 33 aircraft, with solution times within 22.22 to 58.57 minutes—--instances that the earlier model could not solve.
Additionally, in our model, we considered  aircraft that remain in the TMA from the previous time period, a critical aspect for ensuring operational safety that was not addressed by~\cite{spshspp-asmda-21}. 
Furthermore, our framework incorporates both the presence of such aircraft and the requirement of tree consistency.
These results suggest that our approach has the potential to support real-world operations, enabling the regular generation of updated DARs.

%% file: conclusion.tex
\section{Conclusions}\label{con-sec}

We provided a Dantzig-Wolfe reformulation, based on discretisation, of a previous model to compute optimal DARs with fully automated scheduling of CDOs for all aircraft ensuring separation for the arriving traffic.   
While we use CDOs for all aircraft in this paper, the model does not require this, and any length-dependent speed profiles for the aircraft can be used. 
We demonstrated the efficiency of our new model for arrivals in the TMA of Stockholm Arlanda, and for the same grid size as for the previous model: we could significantly improve runtimes 
(from 40.9 hours to within 5 seconds to 12.65 minutes for half-hour high-traffic scenarios) and solve instances that were not solvable with the previous model (33 arriving aircraft in a full hour) within only  22.22 to 58.57 minutes.  
Moreover, the runtimes, in particular, for the 30-minute periods would allow a frequent recomputation of the arrival tree within a reasonable time frame. 
Beyond runtime improvements, our model also includes aircraft carried over from the previous period---a key aspect for operational
 safety not considered by~\cite{spshspp-asmda-21}---as well as constraints ensuring tree consistency.
This clearly shows that our model significantly outperforms the previous model. 

This indicates that in future work, we may be capable to tackle even other influencing factors: 
a fine time resolution and a finer grid resolution.  
The currently used speed profiles only depend on the arrival-path length (determining the distance to go),  and assume no wind and international standard atmospheric conditions. 
This is of course a simplifying assumption, which stems from the previous model's functionality selecting the ``correct'' speed profile based on the path's length, which in turn is also selected in the MIP model.  
This arc- and flow-based approach could not handle selecting the ``correct'' speed profile depending even on the direction of the edges used in the path.  
However,  in our reformulated path-based model,  we select complete paths, and for each aircraft there is a one-to-one correspondence between a path and a speed profile---and we made use of this feature of the reformulation presented in this paper for integrating atmospheric conditions, \cite{hprs-dwraa-25-4}.
Moreover,  for the Stockholm-Arlanda-airport TMA,  with the previous model solutions could be obtained for a $15\times 11$ grid (which---by the edge length---ensures a separation of 6NM), but not for a finer grid.  However, in particular, a finer temporal resolution is of practical interest, because it can reflect actual traversal times more exactly.
For a finer grid or other extensions, we might not be able to include all enumerated path variables, but have to employ column generation.

Additionally, we face uncertainties stemming, for example, from the wind direction (which can be predicted, but not determined with certainty) and actual arrival times of aircraft to the TMA entry points (which can deviate from planned times).  
When we optimize our DARs to specific conditions,  deviations from these could result in infeasible solutions.  Hence, we aim to make the planning even more realistic using robust optimization.

%% file: appendixA.tex
\section{Additional Tables}\label{appA}
\vspace*{-1cm}
\begin{table*}[b!]
    \scriptsize
	\centering
	\caption{Results of experiment set 1a with an upper bound of 15 on the path length}\label{table:results1a-15}
	\begin{tabular}{cccccccc}\hline
                                        &                     &                   &            &                    &  Run      & Trajectory& Average \\
                                        & $|\mathcal{A}|$     &  $|C_2|$          &  $\mu$     & Tree               &  time (s) &generation & deviation \\
                                        &                     &                   &            &                    &           &time (s)   &  \\\hline
        \multirow{14}{*}{5:00--5:29}&\multirow{7}{*}{16} & \multirow{7}{*}{6}& 0          & ---                &  ---      &           &       \\[.2em]
                                        &                    &                    & 1          &$\text{T}^{a}_{33}$ &  143      & 62        &  0.69   \\[.2em]
                                        &                    &                    & 2          &$\text{T}^{a}_{34}$ &  377      & 62        & 1.31    \\[.2em]
                                        &                    &                    & 3          &$\text{T}^{a}_{35}$ &  285      & 64        & 1.87    \\[.2em]
                                        &                    &                    & 4          &$\text{T}^{a}_{36}$ &  1026     & 64        &  2.5   \\[.2em]
                                        &                    &                    & 5          &$\text{T}^{a}_{37}$ &  1162     & 64        &  3.5 \\[.2em]\cmidrule{2-8}
                                        &\multirow{7}{*}{16} & \multirow{7}{*}{0} & 0          &  ---               &   ---     &           &       \\[.2em]
                                        &                    &                    & 1          &$\text{T}^{a}_{38}$ &   110     & 24        & 0.87    \\[.2em]
                                        &                    &                    & 2          &$\text{T}^{a}_{39}$ &  92       & 24        &1.5     \\[.2em]
                                        &                    &                    & 3          &$\text{T}^{a}_{40}$ &  82       & 24        & 2.19    \\[.2em]
                                        &                    &                    & 4          &$\text{T}^{a}_{41}$ &  334      & 24        &2.56     \\[.2em]
                                        &                    &                    & 5          &$\text{T}^{a}_{42}$ &  267      & 24        & 3.62    \\[.2em]\hline
        \multirow{7}{*}{5:30--5:59} &\multirow{4}{*}{17} & \multirow{4}{*}{3} & $\leq 3$   &  ---               &   ---     &           &       \\[.2em]
                                        &                    &                    & 4          &$\text{T}^{a}_{43}$ & 3129      &  77       & 3.06    \\[.2em]
                                        &                    &                    & 5          &$\text{T}^{a}_{44}$ & 2402      &  78       & 3.65   \\[.2em]\cmidrule{2-8}
                                        &\multirow{5}{*}{17} & \multirow{5}{*}{0} & $\leq 2$   &   ---              &  ---      &           &       \\[.2em]
                                        &                    &                    & 3          &$\text{T}^{a}_{45}$ & 3057      & 26        & 2.18    \\[.2em]
                                        &                    &                    & 4          &$\text{T}^{a}_{46}$ & 1186      &  26       & 3    \\[.2em]
                                        &                    &                    & 5          &$\text{T}^{a}_{47}$ & 557       &  26       &  3.88  \\[.2em]\hline
        \multirow{12}{*}{6:00--6:29}&\multirow{6}{*}{9}  & \multirow{6}{*}{2} & $\leq 1$   &   ---              &  ---      &           &       \\[.2em]
                                        &                    &                    & 2          &$\text{T}^{a}_{48}$ & 340       &  36       & 1.44    \\[.2em]
                                        &                    &                    & 3          &$\text{T}^{a}_{49}$ & 511       & 36        & 2.67    \\[.2em]
                                        &                    &                    & 4          &$\text{T}^{a}_{50}$ &  453      & 30        & 3.22    \\[.2em]
                                        &                    &                    & 5          &$\text{T}^{a}_{51}$ & 472       & 34        & 3.67   \\[.2em]\cmidrule{2-8}
                                        &\multirow{6}{*}{9}  & \multirow{6}{*}{0} & $\leq 1$   &   ---              &  --       &           &       \\[.2em]
                                        &                    &                    & 2          &$\text{T}^{a}_{52}$ &  181      & 11        & 1.67    \\[.2em]
                                        &                    &                    & 3          &$\text{T}^{a}_{53}$ &  164      & 11        & 2.67    \\[.2em]
                                        &                    &                    & 4          &$\text{T}^{a}_{54}$ &  195      &11         & 3.22    \\[.2em]
                                        &                    &                    & 5          &$\text{T}^{a}_{55}$ &  215      & 11        &4    \\[.2em]\hline
        \multirow{14}{*}{6:30--6:59}&\multirow{7}{*}{7}  & \multirow{7}{*}{2} & 0          &   ---              &  ---      &           &       \\[.2em]
                                        &                    &                    & 1          &$\text{T}^{a}_{56}$ & 47        & 24        & 1    \\[.2em]
                                        &                    &                    & 2          &$\text{T}^{a}_{57}$ & 79        & 23        & 1.14    \\[.2em]
                                        &                    &                    & 3          &$\text{T}^{a}_{58}$ & 125       & 25        & 1.186    \\[.2em]
                                        &                    &                    & 4          &$\text{T}^{a}_{59}$ & 167       & 26        & 1.71    \\[.2em]
                                        &                    &                    & 5          &$\text{T}^{a}_{60}$ & 102       & 27        & 4   \\[.2em]\cmidrule{2-8}
                                        &\multirow{7}{*}{7} & \multirow{7}{*}{0}  & 0          &  ---               &  ---      &           &       \\[.2em]
                                        &                    &                    & 1          &$\text{T}^{a}_{61}$ & 26        &9          & 0.85    \\[.2em]
                                        &                    &                    & 2          &$\text{T}^{a}_{62}$ & 67        & 9         & 1.28    \\[.2em]
                                        &                    &                    & 3          &$\text{T}^{a}_{63}$ & 77        & 9         & 1.28    \\[.2em]
                                        &                    &                    & 4          &$\text{T}^{a}_{64}$ &  23       & 9         & 3.28    \\[.2em]
                                        &                    &                    & 5          &$\text{T}^{a}_{65}$ &  26       & 9         & 3.57   \\[.2em]\hline
\end{tabular}
\end{table*}

\begin{table*}[b!]
    \scriptsize
	\centering
	\caption{Results of experiment set 1b without light aircraft}	\label{table:results1b-14}
	\begin{tabular}{ccccccccc}\hline
                                       &                   &                     &  Previous                         &           &                    &Run     & Trajectory& Average \\
                                       & $|\mathcal{A}|$   &  $|C_2|$            & tree                              &  $\mu$    &  Tree              &time (s)&generation & deviation \\
                                       &                   &                     &                                   &           &                    &        &time (s)   &  \\\noalign{\hrule height 1pt}
        \multirow{7}{*}{5:30--5:59}&\multirow{7}{*}{17}& \multirow{7}{*}{0}  & $\text{T}^{a}_6$--$\text{T}^{a}_7$& $\leq 8$  &  ---               & ---    &           &     \\[.2em]\cmidrule{4-9}
                                       &                   &                     & $\text{T}^{a}_8$                  & $\leq 4$  & ---                & ---    &           &     \\[.2em]
                                       &                   &                     & $\text{T}^{a}_8$                  & 5         &$\text{T}^{b}_{17}$ & 398    &  9        & 4.59  \\[.2em]\cmidrule{4-9}
                                       &                   &                     & $\text{T}^{a}_9$                  & $\leq 4$  & ---                & ---    &           &     \\[.2em]
                                       &                   &                     & $\text{T}^{a}_9$                  & 5         &$\text{T}^{b}_{18}$ & 377    &  8        & 3.18  \\[.2em]\cmidrule{4-9}
                                       &                   &                     & $\text{T}^{a}_{10}$               & $\leq 4$  & ---                & ---    &           &     \\[.2em]
                                       &                   &                     & $\text{T}^{a}_{10}$               & 5         & $\text{T}^{b}_{19}$& 268    &  8        & 3.06  \\[.2em]\noalign{\hrule height 1pt}
        \multirow{6}{*}{6:00--6:29}&\multirow{6}{*}{9} & \multirow{6}{*}{0}  & $\text{T}^{b}_{19}$               &$\leq 2$   & ---                &  ---   &           &       \\[.2em]
                                       &                   &                     & $\text{T}^{b}_{19}$               & 3         &$\text{T}^{b}_{20}$ & 33     &    3      & 2.44   \\[.2em]
                                       &                   &                     & $\text{T}^{b}_{19}$               & 4         &$\text{T}^{b}_{21}$ & 77     &  3        &  3.33   \\[.2em]
                                       &                   &                     & $\text{T}^{b}_{19}$               & 5         &$\text{T}^{b}_{22}$ &  247   &    4      &  3.78  \\[.2em]
                                       &                   &                     & $\text{T}^{b}_{19}$               & 6         &$\text{T}^{b}_{23}$ &  31    &    4      &  5.78  \\[.2em]\noalign{\hrule height 1pt}
        \multirow{3}{*}{6:30--6:59}&\multirow{3}{*}{7} & \multirow{3}{*}{0} &$\text{T}^{b}_{22}$                 & $\leq 4$  &  ---               &  ---   &           &       \\[.2em]
                                       &                   &                    &$\text{T}^{b}_{22}$                 & 5         &$\text{T}^{b}_{24}$ & 12     &3          & 3.71    \\[.2em]
                                       &                   &                    &$\text{T}^{b}_{22}$                 & 6         &$\text{T}^{b}_{25}$ & 18     & 3         & 4.57    \\[.2em]
                \hline
\end{tabular}
\end{table*}

\begin{table*}[b!]
    \scriptsize
    	\centering
    \caption{Results of experiment set 1c without light aircraft}\label{table:results1c-14}
	\begin{tabular}{cccccccccc}\hline
                                        &                    &                      &   Previous        &          &            &                   &       Run    &  Trajectory & Average \\
                                        & $|\mathcal{A}|$    &  $|C_2|$             & tree              & $\mu$    &  $U$       & Tree              &   time (s)   &generation   & deviation \\
                                        &                    &                      &                   &          &            &                   &              &time (s)     &  \\\noalign{\hrule height 1pt}
        \multirow{30}{*}{5:30--5:59}&\multirow{30}{*}{17}& \multirow{30}{*}{0}  & $\text{T}^{a}_6$  & $\leq 3$ & $\leq 15$  & ---               &   ---        &             &     \\[.2em]
                                        &                    &                      & $\text{T}^{a}_6$  & 4        &$\leq 7$    &---                &  ---         &             &     \\[.2em]
                                        &                    &                      & $\text{T}^{a}_6$  & 4        &8           &$\text{T}^{c}_{18}$& 139          & 8           &  2.94 \\[.2em]
                                        &                    &                      & $\text{T}^{a}_6$  & 5        & $\leq 7$   &---                &  ---         &             &     \\[.2em]
                                        &                    &                      & $\text{T}^{a}_6$  & 5        &   8        &$\text{T}^{c}_{19}$& 44           & 8           & 3.35 \\[.2em]
                                        &                    &                      & $\text{T}^{a}_6$  & 6        & $\leq 4 $  &---                &  ---         &             &     \\[.2em]
                                        &                    &                      & $\text{T}^{a}_6$  & 6        & 5          &$\text{T}^{c}_{20}$&  58          & 8           & 4.88    \\[.2em]\cmidrule{4-10}
                                        &                    &                      & $\text{T}^{a}_7$  & $\leq 3$ & $\leq 15$  & ---               & ---          &             &      \\[.2em]
                                        &                    &                      & $\text{T}^{a}_7$  & 4        & $\leq 7$   &---                &  ---         &             &     \\[.2em]
                                        &                    &                      & $\text{T}^{a}_7$  & 4        & 8          &$\text{T}^{c}_{21}$& 147          & 8           & 2.88   \\[.2em]
                                        &                    &                      & $\text{T}^{a}_7$  & 5        & $\leq 3$   & ---               & ---          &             &     \\[.2em]
                                        &                    &                      & $\text{T}^{a}_7$  & 5        & 4          &$\text{T}^{c}_{22}$& 43           &  8          & 3.76    \\[.2em]
                                        &                    &                      & $\text{T}^{a}_7$  & 6        & 1          &$\text{T}^{c}_{23}$& 11           &  8          & 4.29    \\[.2em]\cmidrule{4-10}
                                        &                    &                      & $\text{T}^{a}_8$  & $\leq 3$ & $\leq 15$  & ---               & ---          &             &      \\[.2em]
                                        &                    &                      & $\text{T}^{a}_8$  & 4        & $\leq 9$   &---                &  ---         &             &     \\[.2em]
                                        &                    &                      & $\text{T}^{a}_8$  & 4        & 10         &$\text{T}^{c}_{24}$&  290         &  8          & 3.06    \\[.2em]
                                        &                    &                      & $\text{T}^{a}_8$  & 5        & 1          &$\text{T}^{c}_{25}$&  9           &   8         & 3.76    \\[.2em]
                                        &                    &                      & $\text{T}^{a}_8$  & 6        & 1          &$\text{T}^{c}_{26}$&  10          &   8         &  4.47   \\[.2em]\cmidrule{4-10}
                                        &                    &                      & $\text{T}^{a}_9$  &$\leq 3$  & $\leq 15$  &---                &  ---         &             &     \\[.2em]
                                        &                    &                      & $\text{T}^{a}_9$  & 4        & $\leq 10$  &---                &  ---         &             &     \\[.2em]
                                        &                    &                      & $\text{T}^{a}_9$  & 4        & 11         &$\text{T}^{c}_{27}$&  289         & 8           & 2.88    \\[.2em]
                                        &                    &                      & $\text{T}^{a}_9$  & 5        & 1          &$\text{T}^{c}_{28}$&  9           & 8           & 3.47    \\[.2em]
                                        &                    &                      & $\text{T}^{a}_9$  & 6        &1           &$\text{T}^{c}_{29}$& 11           & 8           & 4.59 \\[.2em]\cmidrule{4-10}
                                        &                    &                      &$\text{T}^{a}_{10}$&$\leq 3$  & $\leq 15$  &---                &  ---         &             &     \\[.2em]
                                        &                    &                      &$\text{T}^{a}_{10}$& 4        & $\leq 10$  &---                &  ---         &             &     \\[.2em]
                                        &                    &                      &$\text{T}^{a}_{10}$& 4        & 11         &$\text{T}^{c}_{30}$&  289         & 8           & 2.71    \\[.2em]
                                        &                    &                      &$\text{T}^{a}_{10}$& 5        & 1          &$\text{T}^{c}_{31}$&  9           & 8           & 3.65    \\[.2em]
                                                                                     \noalign{\hrule height 1pt}
        \multirow{14}{*}{6:00--6:29}&\multirow{14}{*}{9} & \multirow{14}{*}{0}  &$\text{T}^{c}_{24}$& 5        &  $\leq 9 $ & ---               &  ---         &             &       \\[.2em]
                                        &                    &                      &$\text{T}^{c}_{24}$& 5        &  10        &$\text{T}^{c}_{32}$&  68          &  3          & 2.44     \\[.2em]\cmidrule{4-10}
                                        &                    &                      &$\text{T}^{c}_{25}$& 4        & $\leq 9$   & ---               & ---          &             &       \\[.2em]
                                        &                    &                      &$\text{T}^{c}_{25}$& 4        & 10         &$\text{T}^{c}_{33}$& 38           &  4          & 2.67  \\[.2em]
                                        &                    &                      &$\text{T}^{c}_{25}$& 5        &  $\leq 9$  &  ---              &   ---        &             &     \\[.2em]
                                        &                    &                      &$\text{T}^{c}_{25}$& 5        &   10       &$\text{T}^{c}_{34}$& 67           &   4         &   2.44  \\[.2em]\cmidrule{4-10}
                                        &                    &                      &$\text{T}^{c}_{28}$& 4        &  $\leq 10$ &  ---              &   ---        &             &     \\[.2em]
                                        &                    &                      &$\text{T}^{c}_{28}$& 4        &   11       &$\text{T}^{c}_{35}$& 42           &   4         &   3  \\[.2em]
                                        &                    &                      &$\text{T}^{c}_{28}$& 5        &  $\leq 10$ &  ---              &   ---        &             &     \\[.2em]
                                        &                    &                      &$\text{T}^{c}_{28}$& 5        &   11       &$\text{T}^{c}_{36}$& 58           &   4         &   2.89  \\[.2em]\cmidrule{4-10}
                                        &                    &                      &$\text{T}^{c}_{31}$& 4        &  $\leq 10$ &  ---              &   ---        &             &     \\[.2em]
                                        &                    &                      &$\text{T}^{c}_{31}$& 4        &   11       &$\text{T}^{c}_{37}$& 46           &   4         &  2.78  \\[.2em]
                                        &                    &                      &$\text{T}^{c}_{31}$& 5        &  $\leq 10$ &  ---              &   ---        &             &     \\[.2em]
                                        &                    &                      &$\text{T}^{c}_{31}$& 5        &   11       &$\text{T}^{c}_{38}$& 55           &   4         &   2.89  \\[.2em]
                                        \noalign{\hrule height 1pt}
        \multirow{1}{*}{6:30--6:59} &\multirow{1}{*}{7}  & \multirow{1}{*}{0}   &$\text{T}^{c}_{38}$& 3        & 1          &$\text{T}^{c}_{39}$& 3            &  3          &   1    \\[.2em]
                                       \hline
\end{tabular}
\end{table*}

\begin{table*}[b!]
    \scriptsize
	\centering
	\caption{Results of experiment set 1d without light aircraft}\label{table:results1d-14}
	\begin{tabular}{cccccccccc}\hline
                                        &                     &                         &   Previous          &          &            &                 &       Run    &  Trajectory & Average \\
                                        & $|\mathcal{A}|$     &  $|C_2|$                & tree                & $\mu$    &  $U$       & Tree            &   time (s)   &generation   & deviation \\
                                        &                     &                         &                     &          &            &                 &              &time (s)     &  \\[.2em]\noalign{\hrule height 1pt}
        \multirow{3}{*}{5:30--5:59} &\multirow{3}{*}{17}  & \multirow{3}{*}{0}      & $\text{T}^{a}_{10}$ & $\leq 8$ & $\leq 10$  & ---             &   ---        &             &     \\[.2em]
                                        &                     &                         & $\text{T}^{a}_{10}$ & 9        & $\leq 6$   & ---             &  ---         &             &     \\[.2em]
                                        &                     &                         & $\text{T}^{a}_{10}$ & 9        & 7          &$\text{T}^{d}_4$ &  461         & 8           &  7.12   \\[.2em]
                                \noalign{\hrule height 1pt}
        \multirow{3}{*}{6:00--6:29} &\multirow{3}{*}{9}   & \multirow{3}{*}{0}      &$\text{T}^{d}_4$     &$\leq 9$  & $\leq 15$  & ---             &   ---        &             &     \\[.2em]
                                        &                     &                         &$\text{T}^{d}_4$     & 10       &  $\leq 12$ &  ---            & ---          &             &     \\[.2em]
                                        &                     &                         &$\text{T}^{d}_4$     & 10       &    13      &$\text{T}^{d}_5$ & 112          &   4         &  5.89  \\[.2em]
                            \noalign{\hrule height 1pt}
        \multirow{3}{*}{6:30--6:59} &\multirow{3}{*}{7}   & \multirow{3}{*}{0}      &$\text{T}^{d}_5$     &$\leq 10$ & $\leq 15$  & ---             &   ---        &             &     \\[.2em]
                                        &                     &                         &$\text{T}^{d}_5$     & 11       &  $\leq 6$  &  ---            & ---          &             &     \\[.2em]
                                        &                     &                         &$\text{T}^{d}_5$     & 11       &    7       &$\text{T}^{d}_6$ & 38           &   3         &  8.28  \\[.2em]
                                        \hline
\end{tabular}
\end{table*}

\begin{table*}[b!]
    \scriptsize
	\centering
	\caption{Results of experiment set 2b without light aircraft}\label{table:results2b-14}
	\begin{tabular}{ccccccccc}\hline
                                        &                    &                      &  Previous         &            &                    &     Run     &  Trajectory   &  Average\\
                                        & $|\mathcal{A}|$    &  $|C_2|$             & tree              &  $\mu$     & Tree               &   time (s)  &generation     & deviation \\
                                        &                    &                      &                   &            &                    &             &time (s)       &  \\\hline
        \multirow{10}{*}{6:00--6:59}&\multirow{10}{*}{16}& \multirow{10}{*}{0}  &$\text{T}^{a}_{67}$& $\leq 2$   & ---                &  ---        &               &       \\
                                        &                    &                      &$\text{T}^{a}_{67}$& 3          &$\text{T}^{b}_{29}$ &  40         &   12          &   2.37     \\
                                        &                    &                      &$\text{T}^{a}_{67}$& 4          &$\text{T}^{b}_{30}$ &   55        &   12          &   3      \\
                                        &                    &                      &$\text{T}^{a}_{67}$& 5          &$\text{T}^{b}_{31}$ &   25        &    12         &  4.25      \\
                                        &                    &                      &$\text{T}^{a}_{67}$& 6          &$\text{T}^{b}_{32}$ &   89        &       12      &   5.06     \\\cmidrule{4-9}
                                        &                    &                      &$\text{T}^{a}_{68}$& $\leq 2$   & ---                & ---         &               &        \\
                                        &                    &                      &$\text{T}^{a}_{68}$& 3          &$\text{T}^{b}_{33}$ &  39         &   12          &   2.31     \\
                                        &                    &                      &$\text{T}^{a}_{68}$& 4          &$\text{T}^{b}_{34}$ &   46        &   12          &   3      \\
                                        &                    &                      &$\text{T}^{a}_{68}$& 5          &$\text{T}^{b}_{35}$ &   56        &    12         &  3.87     \\
                                        &                    &                      &$\text{T}^{a}_{68}$& 6          &$\text{T}^{b}_{36}$ &   35        &       12      &   4.94     \\\hline
\end{tabular}
\end{table*}

\begin{table*}[b!]
    \scriptsize
	\centering
	\caption{Results of experiment set 2c without light aircraft}\label{table:results2c-14}
	\begin{tabular}{cccccccccc}\hline
                                        &                    &                      & Previous          &           &           &                    &     Run   & Trajectory&  Average\\
                                        & $|\mathcal{A}|$    &  $|C_2|$             & tree              &  $\mu$    &$U$        & Tree               &   time (s)&generation & deviation \\
                                        &                    &                      &                   &           &           &                    &           &time (s)   &  \\\noalign{\hrule height 1pt}
        \multirow{14}{*}{6:00--6:59}&\multirow{14}{*}{16}& \multirow{14}{*}{0}  &$\text{T}^{a}_{67}$& 3         &$\leq 15$  & ---                &  ---      &           &       \\
                                        &                    &                      &$\text{T}^{a}_{67}$& 4         &$\leq 9$   & ---                &   ---     &           &       \\
                                        &                    &                      &$\text{T}^{a}_{67}$& 4         & 10        &$\text{T}^{c}_{42}$ &   35      &   12      &  3.25     \\
                                        &                    &                      &$\text{T}^{a}_{67}$& 5         &$\leq 9$   & ---                &   ---     &           &       \\
                                        &                    &                      &$\text{T}^{a}_{67}$& 5         & 10        &$\text{T}^{c}_{43}$ & 119        &  12       &  3.44     \\
                                        &                    &                      &$\text{T}^{a}_{67}$& 6         & $\leq 3$  & ---                &  ---      &           &       \\
                                        &                    &                      &$\text{T}^{a}_{67}$& 6         &4          &$\text{T}^{c}_{44}$ & 23        & 12        &   4.87    \\\cmidrule{4-10}
                                        &                    &                      &$\text{T}^{a}_{68}$& 3         &$\leq 15$  & ---                &   ---     &           &       \\
                                        &                    &                      &$\text{T}^{a}_{68}$& 4         &$\leq 10$  & ---                &   ---     &           &       \\
                                        &                    &                      &$\text{T}^{a}_{68}$& 4         & 11        &$\text{T}^{c}_{45}$ &   32      &   12      &  3.25     \\
                                        &                    &                      &$\text{T}^{a}_{68}$& 5         &$\leq 10$  & ---                &   ---     &           &       \\
                                        &                    &                      &$\text{T}^{a}_{68}$& 5         & 11        &$\text{T}^{c}_{46}$ &   32      &   12      &  3.81     \\
                                        &                    &                      &$\text{T}^{a}_{68}$& 6         & $\leq 3$  & ---                &  ---      &           &       \\
                                        &                    &                      &$\text{T}^{a}_{68}$& 6         &4          &$\text{T}^{c}_{47}$ &   31      &   12      &  5.12     \\
                              \hline
\end{tabular}
\end{table*}

\begin{table*}[b!]
    \scriptsize
	\centering
	\caption{Results of experiment set 2d without light aircraft}\label{table:results2d-14}
	\begin{tabular}{cccccccccc}\hline                                        
                                        &                    &                      &  Previous         &           &           &                    &       Run & Trajectory&  Average\\
                                        & $|\mathcal{A}|$    &  $|C_2|$             & tree              &  $\mu$    &$U$        & Tree               &   time (s)&generation & deviation \\
                                        &                    &                      &                   &           &           &                    &           &time (s)   &  \\\hline
        \multirow{14}{*}{6:00--6:59}&\multirow{14}{*}{16}& \multirow{14}{*}{0}  &$\text{T}^{a}_{67}$& 3         &$\leq 15$  & ---                &  ---      &           &       \\
                                        &                    &                      &$\text{T}^{a}_{67}$& 4         &$\leq 13$   & ---                &   ---     &           &       \\
                                        &                    &                      &$\text{T}^{a}_{67}$& 4         & 14        &$\text{T}^{d}_9$    &   37      &   12      &  3.31     \\
                                        &                    &                      &$\text{T}^{a}_{67}$& 5         &$\leq 9$   & ---                &   ---     &           &       \\
                                        &                    &                      &$\text{T}^{a}_{67}$& 5         & 10        &$\text{T}^{d}_{10}$ & 56        &  12       &  3.62     \\
                                        &                    &                      &$\text{T}^{a}_{67}$& 6         & $\leq 3$  & ---                &  ---      &           &       \\
                                        &                    &                      &$\text{T}^{a}_{67}$& 6         &4          &$\text{T}^{d}_{11}$ & 24        & 12        &   5.12    \\\cmidrule{4-10}
                                        &                    &                      &$\text{T}^{a}_{68}$& 3         &$\leq 15$  & ---                &   ---     &           &       \\
                                        &                    &                      &$\text{T}^{a}_{68}$& 4         &$\leq 14$  & ---                &   ---     &           &       \\
                                        &                    &                      &$\text{T}^{a}_{68}$& 4         & 15        &$\text{T}^{d}_{12}$ &   35      &   12      &  3.25     \\
                                        &                    &                      &$\text{T}^{a}_{68}$& 5         &$\leq 10$  & ---                &   ---     &           &       \\
                                        &                    &                      &$\text{T}^{a}_{68}$& 5         & 11        &$\text{T}^{d}_{13}$ &   41      &   12      &  3.81     \\
                                        &                    &                      &$\text{T}^{a}_{68}$& 6         & $\leq 3$  & ---                &  ---      &           &       \\
                                        &                    &                      &$\text{T}^{a}_{68}$& 6         &4          &$\text{T}^{d}_{14}$ &   38      &   12      &  5.06     \\ \hline
\end{tabular}
\end{table*}